\documentclass[preprint,11pt,nopreprintline]{elsarticle}

\usepackage{amsmath}
\usepackage{amsfonts}
\usepackage{amssymb}
\usepackage{amsthm}
\usepackage{bm}

\usepackage{lineno}
\usepackage{multirow}
\usepackage{subfigure}
\usepackage{mathrsfs}
\usepackage[margin=1in]{geometry}
\usepackage[hidelinks]{hyperref}
\usepackage{cleveref}
\crefname{equation}{equation}{equations}

\theoremstyle{plain}
\newtheorem{thm}{Theorem}[section]

\theoremstyle{definition}

\theoremstyle{remark}
\newtheorem{remark}[thm]{Remark}

\let\oldequation\equation
\let\oldendequation\endequation
\renewenvironment{equation}
  {\linenomathNonumbers\oldequation}
  {\oldendequation\endlinenomath}
  
\let\oldalign\align
\let\oldendalign\endalign
\renewenvironment{align}
  {\linenomathNonumbers\oldalign}
  {\oldendalign\endlinenomath}

\begin{document}

\begin{frontmatter}



\title{A Cartesian grid-based boundary integral method for moving interface problems}

\author[inst1]{Han Zhou}\ead{hzhou24@sas.upenn.edu}
\author[inst2]{Shuwang Li}\ead{sli15@illinoistech.edu}
\author[inst3]{Wenjun Ying\corref{corauthor}}\ead{wying@sjtu.edu.cn}

\cortext[corauthor]{Corresponding author}

\address[inst1]{Department of Mathematics, University of Pennsylvania, Philadelphia, PA 19104, USA}
\address[inst2]{Department of Applied Mathematics, Illinois Institute of Technology, Chicago, IL 60616, USA}
\address[inst3]{School of Mathematical Sciences, MOE-LSC and Institute of Natural Sciences, Shanghai Jiao Tong University, Minhang, Shanghai 200240, China}

\begin{abstract}
We develop a Cartesian grid-based boundary integral framework for moving interface problems and apply it to two representative examples: the Hele-Shaw flow and the Stefan problem.
The bulk elliptic and parabolic PDEs are reformulated as boundary integral equations and solved with the matrix-free generalized minimal residual (GMRES) method.
The boundary integrals are evaluated through equivalent interface problems discretized by finite difference methods, which avoids singular and nearly singular quadrature and enables the use of fast PDE solvers such as the fast Fourier transform (FFT) and geometric multigrid methods on fixed Cartesian grids.
To evolve the moving interface, we use the $\theta-L$ formulation together with small-scale decomposition, which maintains mesh quality and removes the stiffness induced by the curvature term through efficient semi-implicit time stepping.
To the best of our knowledge, this is the first Cartesian grid-based numerical framework that successfully combines the small-scale decomposition and spatiotemporal rescaling ideas for moving interface problems.
The resulting framework provides a unified approach for both elliptic- and parabolic-type moving interface problems.
Numerical examples, including long-time Hele-Shaw simulations and dendritic solidification with and without flow, demonstrate the accuracy, stability, and robustness of the method.
\end{abstract}


\end{frontmatter}

\section{Introduction}\label{sec:intro}

Moving interface problems arise in many areas of the natural sciences and engineering, ranging from mathematics \cite{Cao2007,Chen1992,Bellettini1996} and fluid mechanics \cite{HYMAN1984396,Glimm1980a,Glimm1980b} to materials science \cite{Mullins1956,Sethian1992,Meiron1986,Mullins1963} and imaging science \cite{Li2015,Benes2004}. In these problems, interfaces divide the surrounding region into subregions in which the underlying physics is governed by PDEs.

When the interface motion is not known a priori and must be determined as part of the solution, the problem is called a free boundary problem. Free boundary problems are inherently nonlinear because the interface dynamics are coupled to the underlying PDEs.

In this paper, we consider two representative free boundary problems: the Hele-Shaw flow \cite{DEGREGORIA1988,Hou1994,Li2007} and the Stefan problem \cite{Meiron1986,Juric1996,STRAIN1989342}. The Hele-Shaw problem describes the flow of a viscous fluid in a thin gap between two parallel plates and is elliptic in character. The Stefan problem models solidification or melting with a moving interface and is parabolic in character. Both have been studied extensively because of their practical importance and computational difficulty.

Developing accurate and efficient numerical methods for moving interface problems poses several challenges. The first is the accurate representation of a complex, evolving interface. Although many numerical methods have been developed, including front tracking \cite{Zhao2004,Schmidt1996,Juric1996}, level set methods \cite{Chen1997,Kim2000,Gibou2003,Boledi2022,Limare2023}, volume-of-fluid methods \cite{Scardovelli1999,HIRT1981201}, and phase-field methods \cite{Karma1998,Wang2008,Hu2009}, achieving both simplicity and accuracy remains difficult.
The second challenge is the solution of PDEs on complex, time-dependent domains. Because the interface evolves, methods such as the classical finite element method, which rely on body-fitted meshes, require frequent remeshing to maintain mesh quality. This increases both computational cost and implementation complexity. Conventional boundary integral equation methods \cite{STRAIN1989342,Hou1994,Zhu1996,Cristini2002,Cristini2004,Li2007} are efficient for homogeneous PDEs, but they typically require the explicit form of Green's functions for quadrature design together with subtle treatments of singularities.
In recent decades, Cartesian grid-based methods have gained popularity for moving interface problems. These methods, including the immersed boundary method (IBM) \cite{Peskin1977,Peskin2002,Taira2007}, the immersed interface method (IIM) \cite{Li1994,Li1997,Leveque1997,Li2001}, and the ghost fluid method (GFM) \cite{Fedkiw1999,Fedkiw1999a,Liu2000,Nguyen2001}, involve immersing the moving interface into a fixed background mesh, typically a Cartesian grid. This approach simplifies the algorithm and improves computational efficiency.

When considering the effect of surface tension on the interface, the problem formulation includes the Laplace-Young equation or the Gibbs-Thomson relation to account for local curvature. However, the presence of high-order derivatives in the curvature introduces stiffness into the evolution problem and imposes strict stability constraints on the time step when using explicit time-stepping schemes. Conversely, using a straightforward implicit scheme becomes complicated and computationally expensive due to the nonlinear and nonlocal nature of the interface velocity as a function of the interface position.
To address these challenges, Hou et al. developed the small-scale decomposition (SSD) method \cite{Hou1994}. This method removes the stiffness induced by surface tension by combining a special $\theta-L$ formulation of the interface with an implicit discretization of the stiff but linear part of the evolution equation. As a result, it permits much larger time steps.
The SSD method has been adapted for solving multiple moving interface problems, including microstructural evolution in inhomogeneous elastic media \cite{Jou1997}, elastic membranes in viscous flows \cite{Hou2008,Hou2008a}, solid tumor growth \cite{Cristini2003}, and crystal growth \cite{Li2004,Li2005a,Li2005b}, among others.

The goal of this paper is to develop a Cartesian grid-based boundary integral framework for moving interface problems and to apply it to the Hele-Shaw flow and the Stefan problem. After reformulating the bulk PDEs as boundary integral equations, we solve them with the kernel-free boundary integral (KFBI) method, a finite-difference analogue of the classical boundary integral method. The KFBI method is based on potential theory and is designed for elliptic boundary value problems on irregular domains \cite{Ying2007}. It leverages fast PDE solvers on Cartesian grids and avoids the evaluation of singular and nearly singular integrals. Higher-order versions have also been developed and applied successfully to a range of problems \cite{Ying2013,Xie2020,Ying2014,Cao2022}. Compared with traditional boundary integral methods, the present method is quadrature-free for most of the computation and can handle PDEs with variable coefficients \cite{Ying2014,Cao2022}; compared with traditional finite element and finite difference methods on irregular domains, it yields a well-conditioned discrete linear system and fast convergence when GMRES is used. To evolve the moving interface accurately and stably, we use the $\theta-L$ formulation together with the SSD method, which maintains mesh quality and removes curvature-induced stiffness through semi-implicit time stepping. The resulting framework provides a unified treatment of both elliptic- and parabolic-type moving interface problems on fixed Cartesian grids.

The rest of this paper is organized as follows.
Section \ref{sec:prob} introduces the governing equations for the Hele-Shaw and Stefan problems.
Section \ref{sec:BIE} presents the time discretization of the time-dependent PDEs and the corresponding boundary integral formulations.
Section \ref{sec:KFBI} describes the kernel-free boundary integral method.
Section \ref{sec:IEM} presents the numerical treatment of interface evolution.
Section \ref{sec:res} reports numerical results.
Finally, Section \ref{sec:discuss} summarizes the main advantages of the method and discusses possible extensions.
\section{Moving interface problems}\label{sec:prob}
\begin{figure}[htbp]
    \centering
    \includegraphics[width=0.3\textwidth]{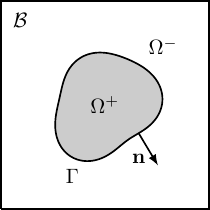}
    \caption{A schematic of a moving interface problem.}
    \label{fig:domains}
\end{figure}

Let $\Gamma:\mathbb R/(2\pi\mathbb Z)\times[0,T]\rightarrow \mathbb{R}^2$ be a time-dependent and closed curve separating the domain of interest $\mathcal{U}\subset \mathbb{R}^2$ into an interior domain $\Omega^+$ and an exterior domain $\Omega^-$. 
The domain $\mathcal{U}$ is either bounded or unbounded depending on the configuration of the problem. 
Throughout this paper, we assume that the interior domain $\Omega^+$ is bounded and contained in a bounding box $\mathcal{B}$; see Fig \ref{fig:domains}.
Physical quantities, such as temperature and velocity fields, satisfy certain PDEs together with interface or boundary conditions imposed on the moving interface $\Gamma$.

\subsection{The Hele-Shaw flow}
The first problem we consider is the Hele-Shaw flow, which describes the motion of the interface between two immiscible viscous fluids in a Hele-Shaw cell—a thin gap between two parallel plates.
The interface, denoted by $\Gamma$, separates the entire space $\mathcal{U}=\mathbb{R}^2$ into an air domain $\Omega^+$ and an oil domain $\Omega^-$.
Let $\mathbf{u}$ and $p$ be the fluid velocity and pressure in the oil domain, respectively. The flow of the oil is assumed to be incompressible and satisfies Darcy's law, namely, 
\begin{equation}\label{eqn:hele-shaw-1}
    \mathbf{u} = -M \nabla p,\quad \nabla\cdot \mathbf{u} = 0, \quad \text{in } \Omega^-,
\end{equation}
where $M=\frac{b^2}{12\mu}$ is the mobility of the fluid, $b$ is the gap width of the Hele-Shaw cell, and $\mu$ is the viscosity.
The pressure in the air domain is assumed to be constant, which we may take to be zero.
External air is injected at the origin, leading to a singular source term:
\begin{equation}\label{eqn:hele-shaw-2}
    \nabla \cdot \mathbf{u} = 2\pi J\delta(\mathbf{x}),\quad \text{in }\Omega^-,
\end{equation}
where $J\geq 0$ is a constant injection rate and $\delta$ is the Dirac delta function.
Combining equations \eqref{eqn:hele-shaw-1} and \eqref{eqn:hele-shaw-2}, and setting $M=1$, we obtain the Poisson equation for the pressure
\begin{equation}\label{eqn:PoissonEqn}
    \Delta p = -2\pi J\delta(\mathbf{x}),\quad \text{in } \Omega^-.
\end{equation}
While the Dirac delta function on the right-hand side of \eqref{eqn:PoissonEqn} vanishes in $\Omega^-$, it prescribes the behavior of the solution at infinity $p = -J\ln |\mathbf{x}| + C + o(1)$ as $|\mathbf{x}|\rightarrow\infty$, where $C$ is a constant ambient pressure.
The surface tension of the air-oil interface is described by the Laplace-Young condition,
\begin{equation} \label{eqn:PoissonBC}
    p =- \sigma \kappa,\quad \text{on }\Gamma,
\end{equation}
where $\sigma>0$ is the surface tension coefficient and $\kappa$ is the curvature of  $\Gamma$.
In addition, the motion of the moving interface $\Gamma$ follows the kinematic condition
\begin{equation}
   \dfrac{d\mathbf{x}}{dt} = \mathbf{u}, \quad \text{for }\mathbf{x}\in\Gamma.
\end{equation}
Since we are only interested in the shape of $\Gamma$, it suffices to prescribe the normal velocity $U$,
\begin{equation}
    U = \frac{d\mathbf{x}}{dt}\cdot\mathbf{n}  = \mathbf{u}\cdot\mathbf{n} = -\partial_{\mathbf{n}}p, \quad \text{for }\mathbf{x}\in\Gamma.
\end{equation}
where $\mathbf{n}$ is the unit outward normal to $\Gamma$.

\subsection{The Stefan problem}
We also consider the Stefan problem, which models the diffusion-driven interface motion between solid and liquid phases. 
The solid-liquid interface $\Gamma$ separates the rectangular domain $\mathcal{B}$ into the solid region $\Omega^+$ and the liquid region $\Omega^-$.
In the classical Stefan problem, thermal convection is ignored in the liquid region, and the temperature field $T$ is assumed to satisfy the heat equation in both the liquid and solid regions,
\begin{align}\label{eqn:heat-eqn}
    \partial_t T = \Delta T,\quad \text{in }\Omega^+\cup\Omega^-,
\end{align}
where we have assumed equal thermal conductivity in both regions.
It is often more realistic to consider convection in the liquid region.
The temperature differences in the fluid lead to changes in the specific volume of fluid parcels and, hence, the fluid density.
The density changes further lead to buoyancy-driven convection of the fluid.
With the Boussinesq approximation \cite{Roger2002}, the fluid is assumed to be incompressible, and the effect of density changes appears only in the buoyancy force.
After incorporating the natural convection effect, the temperature field is described by a coupled system of the (advection-)diffusion equation and the Navier-Stokes equations
\begin{align}
    \partial_t T  = \Delta T, &\quad \text{in }\Omega^+, \\
    \partial_t T + \mathbf{u}\cdot \nabla T = \Delta T, &\quad \text{in }\Omega^-, \label{eqn:AdvDiff}\\
    \partial_t \mathbf{u} + (\mathbf{u}\cdot \nabla)\mathbf{u} = \Delta\mathbf{u}-\nabla p + \mathbf{G}, &\quad\text{in } \Omega^-,\label{eqn:NSE-1}\\
    \nabla \cdot \mathbf{u} = 0, &\quad \text{in } \Omega^-,\label{eqn:NSE-2}
\end{align}
where $\mathbf{u}$ and $p$ are the velocity and pressure of the fluid flow, respectively. Here, we assume that the fluid density and viscosity are unity. The buoyancy force $\mathbf{G}$ under the Boussinesq approximation is linearly proportional to the temperature difference,
\begin{equation}
    \mathbf{G} = -g\beta(T-T_{\infty})\mathbf{j},
\end{equation}
where $\mathbf{j}$ is the unit vector in the vertical direction. The model is a two-way coupling of the temperature and the velocity field through fluid convection and the buoyancy force.

The temperature on the solid-liquid interface is continuous and is prescribed by the Gibbs-Thomson relation,
\begin{align}\label{eqn:GT-rela}
    T +\varepsilon_C(\mathbf{n}) \kappa + \varepsilon_V(\mathbf{n}) U = 0, \quad \text{on } \Gamma, 
\end{align}
where $\varepsilon_C(\mathbf{n})$ and $\varepsilon_V(\mathbf{n})$ are non-negative surface tension and molecular kinetic coefficients. The two coefficients may depend on the orientation of the interface, in which case the problem is anisotropic. 
The thermal flux has a jump across the interface, which is related to the interface velocity through the Stefan condition,
\begin{align}\label{eqn:st-eqn}
    U = [\partial_{\mathbf{n}} T].
\end{align}
Here, the notation $[\cdot]$ denotes the jump of a quantity across the interface.
For example, for a piecewise continuous function $q$, it is defined by
\begin{equation}
    [q](\mathbf{x}) = q|_{\Omega^+}(\mathbf{x}) - q|_{\Omega^-}(\mathbf{x}) = \lim_{\mathbf{y}\in\Omega^+,\mathbf{y}\rightarrow\mathbf{x}}q(\mathbf{y}) - \lim_{\mathbf{y}\in\Omega^-,\mathbf{y}\rightarrow\mathbf{x}}q(\mathbf{y}), \quad \text{for }\mathbf{x}\in\Gamma .  
\end{equation}
A suitable boundary condition should also be given on the outer boundary $\partial\mathcal{B}$, for which we use the no-flux boundary condition $\partial_{\mathbf{n}} T = 0$.
For the fluid flow, the no-slip boundary condition is given on the solid-liquid interface.
On the outer boundary $\partial\mathcal{B}$, we set $\partial_{\mathbf{n}} T = 0$ and $\mathbf{u} = \mathbf{u}_b$ where $\mathbf{u}_b$ is the boundary data describing inflow/outflow or no-slip boundary conditions.

In the modeling of dendritic solidification problems, a solid seed is initially placed in an undercooled surrounding liquid. The temperature of the solid seed is assumed to equal the melt temperature $T_{m}$, and that of the undercooled liquid is given by $T_{\infty}$.
The degree of undercooling is measured by the Stefan number $St = T_{\infty} - T_m$.
\section{Boundary integral equations}\label{sec:BIE}
The PDEs described in the previous section are solved in boundary integral formulations. The Hele-Shaw flow is an elliptic-type moving interface problem, in which the boundary value problem of the Poisson equation can be reformulated as a well-conditioned boundary integral equation. The Stefan problem is a parabolic-type moving interface problem that requires solving time-dependent PDEs, for which we first reduce parabolic PDEs into elliptic PDEs by time discretization and then formulate boundary integral equations for elliptic PDEs. 

\subsection{The Hele-Shaw flow}
In the Hele-Shaw flow, the solution to the Poisson equation \eqref{eqn:PoissonEqn} with the boundary condition \eqref{eqn:PoissonBC} can be represented as the sum of two functions $v$ and $w$, where
\begin{align}
    v(\mathbf{x}) = -J\ln |\mathbf{x}|,
\end{align}
which corresponds to the point-source term and $w(\mathbf{x})$ satisfies an exterior boundary value problem of the Laplace equation
\begin{align}
    \Delta w = 0,&\quad \text{in }\Omega^-,\\
    w =- \sigma\kappa - v, &\quad \text{on }\Gamma. \label{eqn:ex-Lap-bc}
\end{align}
Let $G_0(\mathbf{x}) = (1/2\pi) \ln|\mathbf{x}|$ be the free-space Green's function associated with the Laplacian $\Delta$. The solution $w(\mathbf{x})$ can be represented as a modified double-layer potential
\begin{equation}\label{eqn:w-form}
    w(\mathbf{x}) = (D\varphi)(\mathbf{x}) + \int_{\Gamma}\varphi(\mathbf{y})\,d\mathbf{s}_{\mathbf{y}} = \int_{\Gamma}\varphi(\mathbf{y}) \left(\dfrac{\partial G_0(\mathbf{y}-\mathbf{x})}{\partial\mathbf{n}_{\mathbf{y}}} +1\right) \,d\mathbf{s}_{\mathbf{y}},
\end{equation}
where $\varphi$ is an unknown dipole density function defined on $\Gamma$. 
The boundary integral formulation of $w$ naturally matches the boundary condition at infinity. 
Restricting \eqref{eqn:w-form} on $\Gamma$ and using the boundary condition \eqref{eqn:ex-Lap-bc}, we can obtain a boundary integral equation for the density function $\varphi$,
\begin{equation}\label{eqn:bie-HS}
        -\dfrac{1}{2}\varphi(\mathbf{x}) + \int_{\Gamma}\varphi(\mathbf{y}) \left(\dfrac{\partial G_0(\mathbf{y} -\mathbf{x})}{\partial\mathbf{n}_{\mathbf{y}}} +1\right) \,d\mathbf{s}_{\mathbf{y}} =- \sigma\kappa - v(\mathbf{x}) , \quad \text{for }\mathbf{x}\in\Gamma,
\end{equation}
which is a Fredholm integral equation of the second kind and is well-conditioned.

\subsection{The Stefan problem}
We first discretize the time derivatives in the advection-diffusion equation \eqref{eqn:AdvDiff} and the Navier-Stokes equations \eqref{eqn:NSE-1}, \eqref{eqn:NSE-2}, thereby reducing the problem at each time step to elliptic PDEs. Let $t_n = n\tau$, $n = 0, 1,\cdots, N_T$, be a uniform temporal mesh, where $\tau = T/N_T$ is the time step. For a function $f$, let $f^n$ denote the numerical approximation of $f(t_n)$.

To simplify the numerical approximation of nonlinear advection terms, we start with the semi-Lagrangian formulation,
\begin{align}
    \dfrac{dT}{dt} = \Delta T, &\quad \text{in }\Omega^+\cup\Omega^-, \label{eqn:sl-AdvDif}\\
    \dfrac{d\mathbf{u}}{dt} =  \Delta\mathbf{u}-\nabla p + \mathbf{G}, &\quad\text{in } \Omega^-,\label{eqn:sl-nse-1}\\
    \nabla \cdot \mathbf{u} = 0, &\quad \text{in } \Omega^-,\label{eqn:sl-nse-2}
\end{align}
where $\frac{d}{dt} = \partial_t + \mathbf{u}\cdot \nabla$ is the material derivative. 
Here, we assume that the velocity $\mathbf{u}$ is extended to the solid region with a value of zero, so that the temperature in both domains satisfies the same equation \eqref{eqn:sl-AdvDif}.
A second-order semi-implicit scheme is adopted to discretize \cref{eqn:sl-AdvDif,eqn:sl-nse-1,eqn:sl-nse-2},
\begin{align}
    \dfrac{3T^{n+1}-4\widetilde{T}^n+\widetilde{T}^{n-1}}{2\tau} = \Delta T^{n+1},&\quad \text{in }\Omega^+\cup\Omega^-,\\
    \dfrac{3\mathbf{u}^{n+1}-4\widetilde{\mathbf{u}}^n+\widetilde{\mathbf{u}}^{n-1}}{2\tau} = \Delta \mathbf{u}^{n+1} -\nabla p^{n+1} +2\mathbf{G}^n-\mathbf{G}^{n-1},&\quad\text{in } \Omega^-,\\
    \nabla\cdot\mathbf{u}^{n+1} =0,&\quad\text{in } \Omega^-.
\end{align}
where $\Tilde{\mathbf{u}}^n$, $\Tilde{T}^n$, $\Tilde{\mathbf{u}}^{n-1}$, and $\Tilde{T}^{n-1}$ are the velocities and temperatures at the departure points $\mathbf{x}^n$ and $\mathbf{x}^{n-1}$, respectively.
The scheme treats the buoyancy term explicitly so that the thermal and flow problems are decoupled and can be solved separately.
The departure points $\mathbf{x}^n$ and $\mathbf{x}^{n-1}$ can be found by solving the initial-value problem backward in time,
\begin{align}
    \dfrac{d\mathbf{x}(t)}{dt} = \mathbf{u}(\mathbf{x}(t), t),\quad \mathbf{x}(t_{n+1})= \mathbf{x}_0.
\end{align}
A second-order midpoint method is utilized for computing the positions of the departure points,
\begin{align}
&\mathbf{x}^*=\mathbf{x}_0-\frac{\tau}{2} \mathbf{u}\left(\mathbf{x}_0-\frac{\tau}{2} \mathbf{u}^{n+\frac{1}{2}}, t_{n+\frac{1}{2}}\right), \quad \mathbf{x}^n=\mathbf{x}_0-\tau \mathbf{u}\left(\mathbf{x}^*, t_{n+\frac{1}{2}}\right), \\
&\mathbf{x}^*=\mathbf{x}_0-\tau \mathbf{u}\left(\mathbf{x}_0-\tau \mathbf{u}^n, t_n\right), \quad \mathbf{x}^{n-1}=\mathbf{x}_0-2 \tau \mathbf{u}\left(\mathbf{x}^*, t_n\right).
\end{align}
Off-grid values of $\mathbf{u}$ are computed with cubic Lagrange interpolation and the velocity at $t_{n+\frac{1}{2}}$ is computed with a second-order extrapolation scheme $\mathbf{u}^{n+\frac{1}{2}} = \frac{3}{2}\mathbf{u}^n - \frac{1}{2}\mathbf{u}^{n-1}$. 

After some rearrangement, the thermal problem reduces to the modified Helmholtz equation
\begin{equation}\label{eqn:mH-2}
    (\Delta-\dfrac{3}{2\tau})T^{n+1} = \dfrac{\widetilde{T}^{n-1}-4\widetilde{T}^n}{2\tau},\quad \text{in }\Omega^+\cup\Omega^-,
\end{equation}
subject to interface conditions $[T^{n+1}] = 0$ and $[\partial_{\mathbf{n}}T^{n+1}] = U$ on $\Gamma$ and the outer boundary condition $\partial_{\mathbf{n}} T^{n+1} = 0$ on $\partial\mathcal{B}$.
Let $c = \sqrt{3/(2\tau)}$.
We split $T^{n+1}$ into two parts $T^{n+1} = T_1 +T_2$, where $T_1$ is the solution to the modified Helmholtz equation with an inhomogeneous right-hand side,
\begin{align}
    \Delta T_1 - c^2 T_1 = \dfrac{\widetilde{T}^{n-1}-4\widetilde{T}^n}{2\tau}, \quad &\text{in }\mathcal{B}, \label{eqn:T1-pde}\\
    \partial_{\mathbf{n}} T_1 = 0,\quad &\text{on }\partial\mathcal{B},\label{eqn:T1-bc}
\end{align}
and $T_2$ is the solution to the interface problem with a homogeneous right-hand side,
\begin{align}\label{eqn:st-prob2}
    \Delta T_2 - c^2 T_2 = 0, \quad &\text{in }\Omega^+\cup\Omega^-, \\
    [T_2] = 0, \quad &\text{on }\Gamma,\\
    [\partial_{\mathbf{n}} T_2] = U,\quad &\text{on }\Gamma,\\
    \partial_{\mathbf{n}} T_2 = 0,\quad &\text{on }\partial\mathcal{B}.
\end{align}
Since the temperature field $T$ is continuous across the interface $\Gamma$, the right-hand side of \eqref{eqn:T1-pde} is also continuous.
The function $T_1$ has high regularity, and \cref{eqn:T1-pde,eqn:T1-bc} can therefore be solved accurately with a standard finite difference scheme.
By contrast, $T_2$ is less smooth across $\Gamma$, so standard finite difference schemes generally do not yield satisfactory results. We therefore reformulate this part as a boundary integral equation and solve it with a Cartesian grid-based method.
For each fixed $\mathbf{y}\in\mathcal{B}$, let $G_{c}(\mathbf{x}, \mathbf{y})$ be the Green's function in the box $\mathcal{B}$,
\begin{align}
    (\Delta - c^2)G_{c}(\mathbf{x}, \mathbf{y}) = \delta(\mathbf{x} - \mathbf{y}), \quad &\text{for }\mathbf{x}\in\mathcal{B},\\
    \partial_{\mathbf{n}}G_{c}(\mathbf{x}, \mathbf{y}) = 0,\quad &\text{for }\mathbf{x}\in\partial\mathcal{B}.
\end{align}
The part $T_2$ can be expressed as a single-layer potential
\begin{align}
    T_2(\mathbf{x}) = -(S\psi)(\mathbf{x}) =- \int_{\Gamma} G_{c}(\mathbf{x}, \mathbf{y}) \psi (\mathbf{y}) \,d\mathbf{s}_y, \quad \text{for }\mathbf{x}\in\mathcal{B},
\end{align}
where $\psi = [\partial_{\mathbf{n}} T_2]= U$ is an unknown density function defined on $\Gamma$.
Further, with the Gibbs-Thomson relation \eqref{eqn:GT-rela}, a boundary integral equation can be obtained for the density function $\psi$,
\begin{align}\label{eqn:bie-GT}
    \varepsilon_V \psi (\mathbf{x})  -\int_{\Gamma} G_{c}(\mathbf{x}, \mathbf{y}) \psi (\mathbf{y}) \,d\mathbf{s}_y  = -\varepsilon_C \kappa(\mathbf{x}) - T_1(\mathbf{x}), \quad \text{for }\mathbf{x}\in\Gamma.
\end{align}
For nonzero $\varepsilon_V$, the boundary integral equation \eqref{eqn:bie-GT} is a Fredholm integral equation of the second kind. 
In the case of local thermal equilibrium, namely $\varepsilon_V=0$, the equation degenerates to a first-kind Fredholm integral equation.

The flow problem leads to the modified Stokes equation
\begin{subequations}\label{eqn:mStokes}
\begin{align}
    (\Delta - \dfrac{3}{2\tau})\mathbf{u}^{n+1}- \nabla p^{n+1}  = 
    \dfrac{\widetilde{\mathbf{u}}^{n-1}-4\widetilde{\mathbf{u}}^n}{2\tau} + 2\mathbf{G}^n-\mathbf{G}^{n-1}, &\quad \text{in }\Omega^-, \label{eqn:mStokes-1}\\
    \nabla\cdot\mathbf{u}^{n+1} = 0, &\quad \text{in }\Omega^-,
\end{align}
\end{subequations}
subject to boundary conditions 
\begin{align}
    \mathbf{u}^{n+1} = \mathbf{0} , \quad \text{on }\Gamma, \quad
    \mathbf{u}^{n+1} = \mathbf{u}_b , \quad \text{on }\partial\mathcal{B},
\end{align}
We also split the solution pair into two parts $(\mathbf{u}^{n+1}, p^{n+1}) = (\mathbf{u}_1, p_1)+(\mathbf{u}_2, p_2)$.
The particular solution $(\mathbf{u}_1, p_1)$ satisfies the modified Stokes equation with an inhomogeneous right-hand side
\begin{subequations}
\begin{align}
    (\Delta - c^2)\mathbf{u}_1- \nabla p_1  = \mathbf{f}, &\quad \text{in }\mathcal{B}, \\
    \nabla\cdot\mathbf{u}_1 = 0, &\quad \text{in }\mathcal{B}, \\
    \mathbf{u}_1 = \mathbf{u}_b, &\quad \text{on }\partial\mathcal{B},
\end{align}
\end{subequations}
where $\mathbf{f}$ is given by
\begin{equation}
\mathbf{f} = 
\left\{
    \begin{aligned}
    &\dfrac{\widetilde{\mathbf{u}}^{n-1}-4\widetilde{\mathbf{u}}^n}{2\tau} + 2\mathbf{G}^n-\mathbf{G}^{n-1}, \quad \text{in }\Omega^-,\\
    &2\mathbf{G}^n-\mathbf{G}^{n-1}, \quad \text{in }\Omega^+.
    \end{aligned}
\right.
\end{equation}
Note that $\mathbf{f}$ is a continuous extension of the right-hand side of \eqref{eqn:mStokes-1} due to the no-slip boundary condition on $\Gamma$.
Then the part $(\mathbf{u}_1, p_1)$ has high regularity and a standard finite difference method for the Stokes equation can be applied.
The second part $(\mathbf{u}_2, p_2)$ satisfies an exterior Dirichlet boundary value problem
\begin{align}
    (\Delta - c^2)\mathbf{u}_2- \nabla p_2  = \mathbf{0}, &\quad \text{in }\Omega^-, \\
    \nabla\cdot\mathbf{u}_2 = 0, &\quad \text{in }\Omega^-, \\
    \mathbf{u}_2 = -\mathbf{u}_1, &\quad \text{on }\Gamma,\\
    \mathbf{u}_2 = \mathbf{0}, &\quad \text{on }\partial\mathcal{B}.
\end{align}
For each fixed $\mathbf{x}\in\mathcal{B}$, let $(\mathbf{G}_{\mathbf{u}}(\mathbf{x}, \mathbf{y}), \mathbf{G}_p(\mathbf{x}, \mathbf{y}))$ be the Green's function pair that satisfies 
\begin{align}
    (\Delta - c^2)\mathbf{G}_{\mathbf{u}}(\mathbf{x}, \mathbf{y}) - \nabla \mathbf{G}_p(\mathbf{x}, \mathbf{y}) = \delta(\mathbf{x}-\mathbf{y})\mathbf{I},\quad &\text{for }\mathbf{x}\in\mathcal{B},\\
    \nabla\cdot \mathbf{G}_p(\mathbf{x}, \mathbf{y}) = 0, \quad &\text{for }\mathbf{x}\in\mathcal{B},\\
    \mathbf{G}_{\mathbf{u}}(\mathbf{x}, \mathbf{y}) = \mathbf{0},\quad &\text{for }\mathbf{x}\in\partial\mathcal{B}.
\end{align}
Then the solution pair $(\mathbf{u}_2, p_2)$ can be represented as a double-layer potential $D\boldsymbol{\varphi} = (D_{\mathbf{u}}\boldsymbol{\varphi},D_p\boldsymbol{\varphi})^T$ \cite{Hsiao2021},
\begin{align}
    \mathbf{u}_2 (\mathbf{x}) &= (D_{\mathbf{u}}\boldsymbol{\varphi})(\mathbf{x})  =  \int_{\Gamma} T(\mathbf{G}_{\mathbf{u}}, \mathbf{G}_p)(\mathbf{x}, \mathbf{y}) \boldsymbol{\varphi}(\mathbf{y}) \,d\mathbf{s}_{\mathbf{y}},\label{eqn:u-double-layer}\\
    p_2 (\mathbf{x}) &= (D_p\boldsymbol{\varphi})(\mathbf{x})= 2\int_{\Gamma} \dfrac{\partial \mathbf{G}_p(\mathbf{x}, \mathbf{y})}{\partial \mathbf{n}_{\mathbf{y}}}  \boldsymbol{\varphi}(\mathbf{y}) \,d\mathbf{s}_{\mathbf{y}}.
\end{align}
where $\boldsymbol{\varphi} = (\varphi_1, \varphi_2)^T$ is a vector-valued unknown density function defined on $\Gamma$ and $T(\mathbf{u}, p)=-p\mathbf{n}+(\nabla\mathbf{u}+\nabla\mathbf{u}^T)\mathbf{n}$ is the traction operator.
Restricting \eqref{eqn:u-double-layer} to $\Gamma$ leads to the boundary integral equation
\begin{align}\label{eqn:ms-bie}
    -\dfrac{1}{2}\boldsymbol{\varphi}(\mathbf{x}) + \int_{\Gamma} T(\mathbf{G}_{\mathbf{u}}, \mathbf{G}_p)(\mathbf{x}, \mathbf{y}) \boldsymbol{\varphi}(\mathbf{y}) \,d\mathbf{s}_{\mathbf{y}} = -\mathbf{u}_1(\mathbf{x}), \quad \text{for }\mathbf{x}\in\Gamma.
\end{align}
The boundary integral equation is also a Fredholm integral equation of the second kind and is well-conditioned.

\section{Kernel-free boundary integral method}\label{sec:KFBI}
The boundary integral equations \eqref{eqn:bie-HS}, \eqref{eqn:bie-GT} and \eqref{eqn:ms-bie} are discretized with a potential theory-based Cartesian grid method, the kernel-free boundary integral method.
The main idea of the KFBI method is to use a Cartesian grid-based PDE solver, rather than numerical quadrature, to evaluate the boundary integrals.
Linear algebraic systems resulting from discretizing the boundary integral equations are solved with a Krylov subspace iterative method, the GMRES method \cite{Saad2003_book}.
In each iteration, one only needs to perform a matrix-vector multiplication, which mainly consists of evaluating boundary integrals.
The procedure can be implemented in a matrix-free manner to avoid forming the full matrix.

\subsection{Equivalent interface problems}
Let $\mathcal{A}$ be an elliptic differential operator that can be the Laplacian, modified Helmholtz, or modified Stokes operator.
The evaluation of boundary integral operators associated with $\mathcal{A}$ can be described in the same framework.
According to the classical potential theory, the single-layer potential $S\boldsymbol{\psi}(\mathbf{x})$ and the double-layer potential $D\boldsymbol{\varphi}(\mathbf{x})$ satisfy equivalent interface problems \cite{Hsiao2021,Ying2007}, which can be unified as
\begin{equation}\label{eqn:kfbi-ifp}
\left\{
    \begin{aligned}
    &\mathcal{A} \mathbf{v} = \mathbf{0}, &\quad \text{in }\Omega^+\cup\Omega^-,\\
    &[\boldsymbol{\pi}_D(\mathbf{v})] = \mathbf{\Phi}, &\quad \text{on }\Gamma,\\
    &[\boldsymbol{\pi}_N(\mathbf{v})] = \mathbf{\Psi}, &\quad \text{on }\Gamma,\\
    &\text{boundary conditions}, &\quad \partial\mathcal{B},
    \end{aligned}
\right.
\end{equation}
where $(\boldsymbol{\pi}_D, \boldsymbol{\pi}_N)$ is the Cauchy data pair and is specified in \Cref{tab:cauchy-data}.
\begin{table}[htbp]
\centering
\caption{Cauchy data pair for different elliptic differential operators.}
\label{tab:cauchy-data}
\begin{tabular}{|c|c|c|c|}
\hline
                                                                & $\mathbf{v}$     & $\boldsymbol{\pi}_D(\mathbf{v})$ & $\boldsymbol{\pi}_N(\mathbf{v})$      \\ \hline
$\mathcal{A}=\Delta-c^2$                                     & $u$              & $u$                  & $\partial_{\mathbf{n}} u$ \\ \hline
$\mathcal{A}=((\Delta-c^2)\mathbf{I}-\nabla, \nabla\cdot)^T$ & $(\mathbf{u},p)^T$ & $\mathbf{u}$         & $T(\mathbf{u}, p)$        \\ \hline
\end{tabular}
\end{table}
The functions $\boldsymbol{\Phi}$ and $\boldsymbol{\Psi}$ are given in \Cref{tab:jumps}.
\begin{table}[htbp]
\centering
\caption{Jump relations of the single- and double-layer potentials.}
\label{tab:jumps}
\begin{tabular}{|c|c|c|}
\hline
                     & $\boldsymbol{\Phi}$       & $\boldsymbol{\Psi}$      \\ \hline
$\mathbf{v} =S\boldsymbol{\psi}$ & $\mathbf{0}$           & $-\boldsymbol{\psi}$ \\ \hline
$\mathbf{v} = D\boldsymbol{\varphi}$  & $\boldsymbol{\varphi}$ & $\mathbf{0}$        \\ \hline
\end{tabular}
\end{table}
The boundary condition on $\partial\mathcal{B}$ depends on which Green's function is used.
Since the Hele-Shaw flow is defined on an unbounded domain, the box boundary $\partial\mathcal{B}$ works only as an artificial boundary. 
A natural choice for the boundary condition is to directly use the integral value on $\partial\mathcal{B}$ as a Dirichlet-type boundary condition, namely, $v(\mathbf{x}) = (D\boldsymbol{\varphi})(\mathbf{x})$ for $\mathbf{x}\in \partial\mathcal{B}$, where the boundary integral is approximated with the composite trapezoidal rule.
Due to the periodicity of the integrand function on $\Gamma$, the composite trapezoidal rule is highly accurate.
It is also worth mentioning that no singular integral needs to be handled here since the quadrature points are on $\Gamma$ and target points are on $\partial\mathcal{B}$.
In the Stefan problem, boundary conditions for the temperature or the fluid field depend on the problem at hand and are relatively easier.

The equivalent interface problem \eqref{eqn:kfbi-ifp} is simpler to solve than the original problems since the interface condition can be easily decoupled into jumps of partial derivatives in each direction, which are shown in \ref{sec:app1}.
After solving the interface problem \eqref{eqn:kfbi-ifp} on a Cartesian grid with an efficient PDE solver, one can extract boundary integral values as well as their normal derivatives from the grid function by interpolation.

We stress that the analytical expression of Green's function is only used in the Hele-Shaw flow to convert the infinite boundary condition to a bounded one. 
In other parts, the method remains kernel-free in the sense that the explicit expression of the integral kernel is not required when computing a boundary integral.
Accelerated by fast PDE solvers, the computation of boundary integrals with the KFBI method is rather efficient and is comparable with the fast multipole method, as demonstrated in \cite{Ying2013}.

\subsection{Corrected finite difference scheme}
The interface problem \eqref{eqn:kfbi-ifp} is solved with a finite difference method with additional corrections near the interface.
For simplicity, suppose the computational domain $\mathcal{B}$ is a unit square $(0,1)^2$  and is uniformly partitioned into $N$ intervals in each spatial direction.
The grid nodes are denoted as $(x_i, y_j)$, $x_i = ih$, $y_j=jh$, $i, j = 0,\cdots, N$ where $h = 1/N$ is the mesh size. 

\subsubsection{The modified Helmholtz equation}
Denote by $u_{i,j}$ the numerical approximation of $u(x_i,y_j)$.
A standard second-order five-point central difference scheme for the modified Helmholtz equation is given by
\begin{align}\label{eqn:fp-fdm}
\Delta_h u_{i,j}  - c^2 u_{i,j}= \dfrac{u_{i+1,j} + u_{i-1,j} + u_{i,j+1} + u_{i,j-1} - 4u_{i,j}}{h^2} - c^2 u_{i,j} = 0.
\end{align}
A grid node $(x_i, y_j)$ is called an irregular node if at least one of the stencil nodes is on the other side of the interface.
Clearly, irregular nodes are always in the vicinity of $\Gamma$.
Since the solution has certain jumps across the interface $\Gamma$, the local truncation error at an irregular node is of order $\mathcal{O}(h^{-2})$ and is too large to achieve an accurate approximation.
To show this, we assume that at an irregular node $(x_i,y_j)\in\Omega^+$, a stencil node $(x_{i+1},y_j)\in\Omega^-$ is on the other side of $\Gamma$ and that the other stencil nodes are all in $\Omega^+$.
Suppose the interface $\Gamma$ intersects the grid line segment between $(x_i,y_j)$ and $(x_{i+1},y_j)$ at $(\xi, y_j)$.
By Taylor expansion, the local truncation error at $(x_i,y_j)$ is given by
\begin{align}\label{eqn:fdm-LTE}
E_h(x_i,y_j) = (\Delta_h -c^2) u(x_i,y_j) = - \dfrac{1}{h^2}\{[u] + (x_{i+1}-\xi)[u_x] + \dfrac{1}{2}(x_{i+1}-\xi)^2[u_{xx}])\} + \mathcal{O}(h).
\end{align}
To avoid the large local truncation error, one can adopt the leading term as a correction term to the right-hand side of \eqref{eqn:fp-fdm}.
This gives the corrected scheme
\begin{align}\label{eqn:modified-fdm}
    \Delta_h u_{i,j} - c^2 u_{i,j} = C_{i,j} = - \dfrac{1}{h^2}\{[u] + (x_{i+1}-\xi)[u_x] + \dfrac{1}{2}(x_{i+1}-\xi)^2[u_{xx}])\}.
\end{align}
The resulting local truncation error becomes $\mathcal{O}(h)$.
Here, the correction term $C_{i,j}$ comes from the contribution of the intersection point $(\xi, y_j)$ and is a linear combination of the jump values $[u]$, $[u_x]$, and $[u_{yy}]$, which are known before solving the interface problem; see \ref{sec:app1}.
Similarly, one can derive correction terms for different intersection patterns by Taylor expansions.
Evidently, $C_{i,j}$ is non-zero only at irregular nodes and is always a linear combination of the jump values $[u],[u_x],[u_y],\cdots$.

Since irregular nodes are near the interface $\Gamma$, which is a co-dimension one object, the overall accuracy can still be second-order \cite{Beale2006}.
Since correction terms only appear on the right-hand side, the coefficient matrix is not altered and FFT-based fast Poisson solvers can be applied to solve the resulting linear system efficiently.

\subsubsection{The modified Stokes equation}
The discretization of the modified Stokes equations is based on the marker and cell (MAC) scheme on a staggered grid.
The pressure is at the cell center, the x-component of the velocity is at the center of the east and west edges of a cell, and the y-component of the velocity is at the center of the north and south edges of a cell.
The discrete mesh functions $u_{i,j}, v_{i,j}, p_{i,j}$ are defined as
\begin{align}
    &u_{i,j} = u(x_i, y_j - h/2), \quad i=0,1\cdots, N_x, \quad j = 1,2\cdots, N_y,\\
    &v_{i,j} = v(x_i - h/2, y_j), \quad i=1,2\cdots, N_x, \quad j = 0,1\cdots, N_y,\\
    &p_{i,j} = p(x_i - h/2, y_j - h/2), \quad i= 1,2\cdots, N_x, \quad j = 1,2\cdots, N_y.
\end{align}
The MAC scheme for the modified Stokes equation is given by
\begin{subequations}\label{eqn:mac-fds}
\begin{align}
    (\Delta_h - c^2) \mathbf{u}_{i,j} - \nabla_h p_{i,j} &= \mathbf{0}, \\
    \nabla_h \cdot \mathbf{u}_{i,j} &= 0,
\end{align}
\end{subequations}
where $\Delta_h, \nabla_h$ are the discrete approximations for $\Delta, \nabla$, respectively, 
\begin{align}
    &\Delta_h \mathbf{u}_{ij} = \dfrac{\mathbf{u}_{i+1,j} + \mathbf{u}_{i-1,j} + \mathbf{u}_{i,j+1} + \mathbf{u}_{i,j-1} -4\mathbf{u}_{i,j} }{h^2},\\
    &\nabla_h p_{i,j} = (\dfrac{p_{i+1,j} -p_{i,j}}{h},\dfrac{p_{i,j+1} -p_{i,j}}{h})^T, \\
    &\nabla_h\cdot \mathbf{u}_{i,j} = \dfrac{u_{i,j}-u_{i-1,j}}{h} + \dfrac{v_{i,j}-v_{i,j-1}}{h}.
\end{align}
The Dirichlet boundary condition is discretized with a symmetric approach for the $x$-component velocity on the top and bottom boundary and the $y$-component velocity on the right and left boundary.
In the presence of interfaces, one can also derive local truncation errors at irregular nodes for each finite difference equation in \eqref{eqn:mac-fds}.
Similarly, by using the leading term in the local truncation error as a correction term of the right-hand side in \eqref{eqn:mac-fds}, we arrive at the corrected MAC scheme
\begin{subequations}\label{eqn:crc-mac-fds}
\begin{align}
    (\Delta_h - c^2) \mathbf{u}_{i,j} - \nabla_h p_{i,j} &= \mathbf{C}_{i,j}, \\
    \nabla_h \cdot \mathbf{u} &= D_{i,j},
\end{align}
\end{subequations}
where $\mathbf{C}_{i,j}$ and $D_{i,j}$ are correction terms, which are linear combinations of jump values $[u]$, $[v]$, $[p]$, $[u_x]$, $[u_y]$, $\cdots$.
The resulting linear system \eqref{eqn:crc-mac-fds} is solved with an efficient V-cycle geometric multigrid method with the Distributive Gauss-Seidel (DGS) smoother \cite{BRANDT197953,BRANDT1992151}.

\begin{remark}
Since particular solutions, such as $T_1$, $(\mathbf{u}_1,p_1)$, have high regularity, classical numerical approaches, such as finite difference methods, finite element methods, and spectral methods, can be applied without corrections.
In this work, we apply the five-point central difference scheme to solve $T_1$ and the MAC scheme to solve $(\mathbf{u}_1,p_1)$. 
\end{remark}

\subsection{Extracting boundary values from the grid data }
The single- and double-layer potential functions $S\boldsymbol{\psi}$ and $D\boldsymbol{\varphi}$ defined by boundary integrals are not smooth on $\Gamma$.
Define the single- and double-layer boundary integral operators $\mathcal{S}\boldsymbol{\psi}$ and $\mathcal{D}\boldsymbol{\varphi}$ as the values of the potential functions $S\boldsymbol{\psi}$ and $D\boldsymbol{\varphi}$ on $\Gamma$, respectively.
Assuming that $\Gamma$ is sufficiently smooth, potential theory implies that the integral operators satisfy the following relations:
\begin{align}
\mathcal{S}\boldsymbol{\psi}(\mathbf{x})& = \dfrac{1}{2}((S\boldsymbol{\psi})^+(\mathbf{x})+ (S\boldsymbol{\psi})^-(\mathbf{x})), \quad \mathbf{x}\in\Gamma,\\
\mathcal{D}\boldsymbol{\varphi}(\mathbf{x}) &= \dfrac{1}{2} ((D\boldsymbol{\varphi})^+(\mathbf{x}) + (D\boldsymbol{\varphi})^-(\mathbf{x})), \quad \mathbf{x}\in\Gamma.
\end{align}
The normal derivatives of the potential functions $\partial_{\mathbf{n}}(S\boldsymbol{\psi})$ and $\partial_{\mathbf{n}}(D\boldsymbol{\varphi})$ have similar expressions:
\begin{align}
    \partial_{\mathbf{n}}(S\boldsymbol{\psi})(\mathbf{x})& = \mathbf{n}(\mathbf{x})\cdot\dfrac{1}{2}((\nabla S\boldsymbol{\psi})^+(\mathbf{x})+ (\nabla S\boldsymbol{\psi})^-(\mathbf{x})), \quad \mathbf{x}\in\Gamma,\\
\partial_{\mathbf{n}}(D\boldsymbol{\varphi}) (\mathbf{x})&= \mathbf{n}(\mathbf{x})\cdot\dfrac{1}{2}((\nabla D\boldsymbol{\varphi})^+(\mathbf{x}) + (\nabla D\boldsymbol{\varphi})^-(\mathbf{x})), \quad \mathbf{x}\in\Gamma.
\end{align}
The normal derivatives of the potential functions are related to the adjoint double-layer integral and hypersingular integral operators \cite{Hsiao2021}.
Once the grid-valued single- and double-layer potential functions $(S\boldsymbol{\psi})_h$ and $(D\boldsymbol{\varphi})_h$ are obtained, boundary integral operators can be computed by extracting boundary values and normal derivatives of the potential functions from grid data via polynomial interpolation.

When performing the interpolation, one needs to be aware of the jump values of the potential functions.
Suppose we evaluate the one-sided limit of a discontinuous function $v(\mathbf{x})$ at a point $\mathbf{p} = (p_1,p_2)^T\in\Gamma$ from the side of $\Omega^+$.
Denote by $v^+(\mathbf{x})$ and $v^-(\mathbf{x})$ two functions that coincide with $v(\mathbf{x})$ in $\Omega^+$ and $\Omega^-$, respectively.
Given a set of interpolation points $\{\mathbf{q}_m\}_{m=1}^M$, by Taylor expansion at $\mathbf{p}$, we have
\begin{align}
v(\mathbf{q}_m)  &= v^+(\mathbf{p}) + (\Delta\mathbf{x}_m)^T\nabla v^+(\mathbf{p}) + \dfrac{1}{2} (\Delta\mathbf{x}_m)^T\nabla^2 v^+ (\mathbf{p})\Delta\mathbf{x}_m + \mathcal{O}(|\Delta\mathbf{x}_m|^3),\quad \text{ if }\mathbf{q}_m\in\Omega^+,\\
v(\mathbf{q}_m) + C(\mathbf{q}_m) &= v^+(\mathbf{p}) + (\Delta\mathbf{x}_m)^T\nabla v^+(\mathbf{p}) + \dfrac{1}{2} (\Delta\mathbf{x}_m)^T\nabla^2 v^+ (\mathbf{p})\Delta\mathbf{x}_m + \mathcal{O}(|\Delta\mathbf{x}_m|^3),\quad \text{ if }\mathbf{q}_m\in\Omega^-,
\end{align}
where $\Delta \mathbf{x}_m = \mathbf{q}_m - \mathbf{p}, m = 1,\cdots,M$ and $C(\mathbf{p})$ is given by
\begin{align}
C(\mathbf{q}_m) = [v](\mathbf{p}) + (\Delta\mathbf{x}_m)^T[\nabla v](\mathbf{p}) + \dfrac{1}{2} (\Delta\mathbf{x}_m)^T[\nabla^2 v] (\mathbf{p})\Delta\mathbf{x}_m.
\end{align}
Dropping the $\mathcal{O}(|\Delta \mathbf{x}_m|^3)$ term, one obtains an $M\times M$ linear system to be solved for the boundary values $v^+$, $\nabla v^+$, and $\nabla^2 v^+$.
Then, the values $v^-$, $\nabla v^-$, and $\nabla^2 v^-$ can also be obtained after simple algebraic manipulations.
Note that a set of interpolation points needs to be selected carefully to ensure the invertibility of the interpolation matrix. To find interpolation points, for example, one first finds the grid node $(x_{i_0}, y_{j_0})$ closest to $\mathbf{p}$.
Let $d_1, d_2$ be two integers,
\begin{equation}
d_1 = 
\left\{
\begin{aligned}
1, \text{ if } p_1 > x_{i_0},\\
-1, \text{ if } p_1 \leq x_{i_0},\\
\end{aligned}
\right.\quad
d_2 = 
\left\{
\begin{aligned}
1, \text{ if } p_2 > y_{j_0},\\
-1, \text{ if } p_2\leq y_{j_0},\\
\end{aligned}
\right.
\end{equation}
Then six interpolation points are chosen
\begin{equation}
\{(x_{i_0+r}, y_{j_0 + s})\}_{(r,s)\in\mathcal{I}},\quad \mathcal{I} = \{(0,0),(1,0),(0,1),(0,-1),(-1,0),(d_1,d_2)\}.
\end{equation}
This choice of interpolation points results in a non-singular $6\times 6$ system, which can be solved with a direct method, such as the LU decomposition.

\section{Interface evolution method}\label{sec:IEM}
\subsection{The $\theta - L$ approach}
Suppose the interface $\Gamma$ is parameterized by $\mathbf{X}(\alpha,t) = (x(\alpha,t),y(\alpha,t))$, where $\alpha\in \mathbb R/(2\pi\mathbb Z)$ is the parameter. Denote by $s = s(\alpha, t)$ the arc-length parameter of the curve.
The evolution of the curve is described by
\begin{align}\label{eqn:ori-evo}
    \partial_t\mathbf{X} = U\mathbf{n} + V\mathbf{s},
\end{align}
where $\mathbf{n}$ is the unit outward normal, $\mathbf{s}$ is the unit tangent, $U$ and $V$ are the normal and tangential components of the curve velocity, respectively. In the Hele-Shaw flow and the Stefan problem, only the shape of the moving interface is of practical interest; its shape can be determined solely by the normal motion, as the tangential motion only reparameterizes the curve.
However, pure normal motion for tracking a moving interface is not a good choice for numerical computation due to the clustering or spread of marker points on the interface.
In the meantime, curvature may introduce numerical stiffness and cause severe stability constraints on the time step.
Instead of using the original $x-y$ description, we adopt the $\theta-L$ description for the interface evolution, that is,
\begin{align}
    L_t &= \int_0^{2\pi} \theta_{\alpha^{\prime}} U \,d\alpha^{\prime}, \label{eqn:L-evo}\\
    \theta_t &= \left(\dfrac{2\pi}{L}\right) (\theta_{\alpha}V - U_{\alpha}),\label{eqn:the-evo}
\end{align}
where $\theta$ is the tangent angle to the curve $\Gamma$ and $L$ is the curve length.
By setting $s_{\alpha}=\sqrt{x_{\alpha}^2+y_{\alpha}^2}\equiv L/2\pi$, one obtains the equal-arclength tangential velocity
\begin{equation}
    V(\alpha,t) = \dfrac{\alpha}{2\pi} \int_{0}^{2\pi} \theta_{\alpha^{\prime}} U \,d\alpha^{\prime} - \int_0^{\alpha} \theta_{\alpha^{\prime}} U \,d\alpha^{\prime}.
\end{equation}
The mapping from $(\theta,L)$ to $(x,y)$ still needs two more integration constants to determine the position of the curve, for which we track $ \mathbf{\overline{X}} = \frac{1}{2\pi}\int_{0}^{2\pi}\mathbf{X}\,d\alpha$ using the evolution equation
\begin{align}\label{eqn:ave-evo}
    \mathbf{\overline{X}}_t = \dfrac{1}{2\pi}\int_{0}^{2\pi} U\mathbf{n} + V\mathbf{s}\,d\alpha.
\end{align}
While it is possible to track a single point, as demonstrated in \cite{Hou1994}, we find that using the averaged value performs better in preserving the symmetry of the interface shape.

\subsection{Small-scale decomposition}
\subsubsection{Hele-Shaw flow}
In the Hele-Shaw flow, the mapping from curvature $\kappa$ to the normal velocity $U$ is a one-phase Dirichlet-to-Neumann (DtN) mapping \cite{pruss2016moving}. 
Similar to \cite{Hou1994}, we apply a small-scale decomposition for the normal velocity.
For a fixed point $\mathbf{x}\in\Gamma$, by Green's third identity, we have
\begin{equation}
    -\int_{\Gamma} \dfrac{\partial}{\partial \mathbf{n}^{\prime}} G_0(\mathbf{x} - \mathbf{x}^{\prime}) (p|_{\Gamma^-})(\mathbf{x}^{\prime}) \,d s_{\mathbf{x}^{\prime}} + \int_{\Gamma} G_0(\mathbf{x} - \mathbf{x}^{\prime})(\partial_{\mathbf{n}}p|_{\Gamma^-})(\mathbf{x}^{\prime}) \,d s_{\mathbf{x}^{\prime}} = \dfrac{1}{2} p|_{\Gamma^-}(\mathbf{x}),
\end{equation}
where $\cdot|_{\Gamma^-}$ denotes the one-sided limit on $\Gamma$ from $\Omega^-$. Since $p|_{\Gamma^-}=-\sigma \kappa = -\sigma \theta_{s}$ and $\partial_{\mathbf{n}}p|_{\Gamma^-}=-U$, we have
\begin{equation}
    \int_{\Gamma} G_0(\mathbf{x} - \mathbf{x}^{\prime})U(\mathbf{x}^{\prime}) \,d s_{\mathbf{x}^{\prime}} = \dfrac{1}{2} \sigma\theta_{s} + \int_{\Gamma} \dfrac{\partial}{\partial \mathbf{n}^{\prime}} G_0(\mathbf{x} - \mathbf{x}^{\prime}) (\sigma\theta_{s})(\mathbf{x}^{\prime}) \,d s_{\mathbf{x}^{\prime}} ,\quad \mathbf{x}\in\Gamma.
\end{equation}
Since the double-layer integral operator $\mathcal{D}$ is a smoothing 
operator \cite{Hsiao2021}, for any function $f\in C(\Gamma)$, we have $(f + \mathcal{D} f) \sim f$. 
Here, the notation $\sim$ means if $f \sim g$ then $f-g$ is smoother than  $f$ and $g$.
Therefore, we have
\begin{equation}\label{hs-eqn:sin-equiv}
    \int_{\Gamma} G_0(\mathbf{x} - \mathbf{x}^{\prime})U(\mathbf{x}^{\prime}) \,d s_{\mathbf{x}^{\prime}} \sim \dfrac{1}{2}\sigma\theta_{s}.
\end{equation}
Taking tangential derivatives on both sides of \eqref{hs-eqn:sin-equiv}, one gets
\begin{align}
\begin{aligned}
    \partial_{s}\int_{\Gamma} G_0(\mathbf{x} - \mathbf{x}^{\prime})U(\mathbf{x}^{\prime}) \,d s_{\mathbf{x}^{\prime}} &= \dfrac{1}{2\pi}\text{p.v.}\int_{0}^{2\pi} \dfrac{(\mathbf{X}(\alpha)-\mathbf{X}(\alpha^{\prime}))\cdot \partial_{\alpha}\mathbf{X}(\alpha)}{|\mathbf{X}(\alpha)-\mathbf{X}(\alpha^{\prime})|^2} U(\mathbf{X}(\alpha^{\prime})) \dfrac{s_{\alpha}(\alpha^{\prime})}{s_{\alpha}(\alpha)} \,d \alpha^{\prime}  \\
    &:= \dfrac{1}{2\pi}\text{p.v.}\int_{0}^{2\pi} I(\alpha,\alpha^{\prime}) \,d\alpha^{\prime}.
\end{aligned}
\end{align}
The integrand $I(\alpha,\alpha^{\prime})$ is singular at  $\alpha^{\prime} = \alpha$. By Taylor expansion, $\mathbf{X}(\alpha)-\mathbf{X}(\alpha^{\prime}) = (\alpha-\alpha^{\prime})\partial_{\alpha}\mathbf{X}(\alpha) + \mathcal{O}(\Delta \alpha^2)$ where $\Delta \alpha=\alpha-\alpha^{\prime}$. Then, we have
\begin{equation}
    I(\alpha,\alpha^{\prime}) \approx \dfrac{1}{\alpha-\alpha^{\prime}} U(\mathbf{X}(\alpha^{\prime})) \approx \dfrac{1}{2}\cot\left(\dfrac{\alpha-\alpha^{\prime}}{2}\right)U(\mathbf{X}(\alpha^{\prime})), \quad \text{as }\alpha^{\prime} \rightarrow \alpha. 
\end{equation}
Therefore, one gets
\begin{align}
    \partial_{s}\int_{\Gamma} G_0(\mathbf{x} - \mathbf{x}^{\prime})U(\mathbf{x}^{\prime}) \,d s_{\mathbf{x}^{\prime}} \sim \dfrac{1}{4\pi}\text{p.v.}\int_{0}^{2\pi}\cot\left(\dfrac{\alpha-\alpha^{\prime}}{2}\right)U(\mathbf{X}(\alpha^{\prime})) \, d\alpha^{\prime} = \dfrac{1}{2}\mathscr{H}[U],
\end{align}
where $\mathscr{H}$  is the Hilbert transform, defined by
\begin{equation}
    \mathscr{H}[\omega] (\alpha) = \dfrac{1}{2\pi} \text{p.v.}\int_{0}^{2\pi}\cot\left(\dfrac{\alpha-\alpha^{\prime}}{2}\right)\omega(\alpha^{\prime}) \, d\alpha^{\prime}.
\end{equation}
Thus, we have
\begin{align}
    \mathscr{H}[U] \sim \sigma\theta_{ss}.
\end{align}
Applying the property of the Hilbert transform $\mathscr{H}[\mathscr{H}[\omega]] = -\omega$, we arrive at
\begin{align}
    U \sim -\sigma \mathscr{H}[\theta_{ss}] = -\sigma \left(\dfrac{2\pi}{L}\right)^2 \mathscr{H}[\theta_{\alpha\alpha}].
\end{align}

\subsubsection{Stefan problem}
For the Stefan problem, the mapping $\varepsilon_C \kappa\rightarrow U$ is a two-phase DtN mapping.
Thus, we derive the small-scale decomposition of $U$.
For simplicity, we assume $\varepsilon_V$ and $\varepsilon_C$ are constants.
Due to the relation $U = \psi$, the DtN mapping is defined through the boundary integral equation \eqref{eqn:bie-GT}. 
Since the integral equation \eqref{eqn:bie-GT} becomes a Fredholm integral equation of the first kind if $\varepsilon_V = 0$, we need to consider two cases: (1) $\varepsilon_V = \mathcal{O}(1)$ and (2) $\varepsilon_V\ll 1$.

\paragraph{Case I: $\varepsilon_V = \mathcal{O}(1)$}
Denote by $\mathcal{S}\psi = (1/\varepsilon_V) \int_{\Gamma}G_c(\mathbf{x}, \mathbf{y}) \psi (\mathbf{y})\, d\mathbf{s}_{\mathbf{y}}$ and $g = -(\varepsilon_C/\varepsilon_V) \kappa - T_1/\varepsilon_V$.
We rewrite the boundary integral equation \eqref{eqn:bie-GT} as an operator equation
\begin{equation}
    (\mathcal{I} - \mathcal{S}) \psi = g,
\end{equation}
where $\mathcal{S}$ is a compact operator and $\mathcal{I} - \mathcal{S}$ has a bounded inverse. In fact, the integral operator $\mathcal{S}$ is a pseudo-differential operator of order $-1$. For a constant $\beta$, it has the mapping property $\mathcal{S}:C^{\beta}(\Gamma)\rightarrow C^{1+\beta}(\Gamma)$. Hence, we have
\begin{align}
    \psi = g + \mathcal{S}\psi = g + S[(\mathcal{I} - \mathcal{S})^{-1}g].
\end{align}
Note that $\mathcal{S}[(\mathcal{I} - \mathcal{S})^{-1}g]$ is smoother than $g$, and $T_1$ is equivalent to a volume integral and is smoother than $\kappa = \theta_s = (2\pi/L)\theta_{\alpha}$. Thus,
\begin{equation}\label{eqn:st-ssd1}
U  \sim - \dfrac{\varepsilon_C}{\varepsilon_V} \left(\dfrac{2\pi}{L}\right) \theta_{\alpha}.
\end{equation}

\paragraph{Case II: $\varepsilon_V \ll 1$}
Define the difference of the two Green's functions $\widetilde{G}_c = G_c (\mathbf{x}, \mathbf{y})- G_0(\mathbf{x}-\mathbf{y})$. Note that, for any $\mathbf{x}\in\mathcal{B}$, the difference satisfies the following problems
\begin{align}
    \Delta \widetilde{G}_{c}(\mathbf{x}, \mathbf{y}) &= c^2 G_c(\mathbf{x}, \mathbf{y}), \quad \mathbf{y}\in\mathcal{B},\\
    \partial_{\mathbf{n}_{\mathbf{y}}}\widetilde{G}_c(\mathbf{x}, \mathbf{y}) &= -\partial_{\mathbf{n}_{\mathbf{y}}}G(\mathbf{x}, \mathbf{y}),\quad \mathbf{y}\in\partial\mathcal{B},
\end{align}
and
\begin{align}
    (\Delta-c^2) \widetilde{G}_{c}(\mathbf{x}, \mathbf{y}) &= c^2 G(\mathbf{x}, \mathbf{y}), \quad \mathbf{y}\in\mathcal{B},\\
    \partial_{\mathbf{n}_{\mathbf{y}}}\widetilde{G}_c(\mathbf{x}, \mathbf{y}) &= -\partial_{\mathbf{n}_{\mathbf{y}}}G(\mathbf{x}, \mathbf{y}),\quad \mathbf{y}\in\partial\mathcal{B}.
\end{align}
By the regularity theory of elliptic PDEs, $\widetilde{G}_c$ has higher regularity than $G_c$ and $G$. Thus, we have
\begin{equation}
    \int_{\Gamma} G_c(\mathbf{x}, \mathbf{y}) \psi(\mathbf{y}) \,d\mathbf{S}_y = \int_{\Gamma} G(\mathbf{x}, \mathbf{y}) \psi(\mathbf{y}) \,d\mathbf{S}_y + \int_{\Gamma} \widetilde{G}_c(\mathbf{x}, \mathbf{y}) \psi(\mathbf{y}) \,d\mathbf{S}_y \sim \int_{\Gamma} G(\mathbf{x}, \mathbf{y}) \psi(\mathbf{y}) \,d\mathbf{S}_y .
\end{equation}
On noting that $\kappa=\theta_s$, from the boundary integral equation \eqref{eqn:bie-GT}, we have
\begin{equation}
    \int_{\Gamma} G(\mathbf{x}, \mathbf{y}) U(\mathbf{y}) \,d\mathbf{S}_y  \sim \varepsilon_C \theta_s.
\end{equation}
Similar to the case of Hele-Shaw flow, taking tangential derivatives on both sides of the equation leads to
\begin{equation}
    \partial_s\int_{\Gamma} G(\mathbf{x}, \mathbf{y}) U(\mathbf{y}) \,d\mathbf{S}_y \sim \dfrac{1}{2} \mathscr{H}[U].
\end{equation}
Recalling that $\mathscr{H}[\mathscr{H}\omega] = -\omega$, we obtain
\begin{equation}\label{eqn:st-ssd2}
    U \sim -2 \mathscr{H}[\partial_s(\varepsilon_C\partial_s\theta)] = -2\varepsilon_C\left(\dfrac{2\pi}{L}\right)^2 \mathscr{H}[\theta_{\alpha\alpha}].
\end{equation}
Thus, we obtain the small-scale decomposition of $U$:
\begin{equation}
    U \sim 
    \left\{
    \begin{aligned}
     &- \dfrac{\varepsilon_C}{\varepsilon_V} \left(\dfrac{2\pi}{L}\right) \theta_{\alpha}, &\quad \text{if } \varepsilon_V=\mathcal{O}(1), \\
    & -2\varepsilon_C\left(\dfrac{2\pi}{L}\right)^2 \mathscr{H}[\theta_{\alpha\alpha}], &\quad \text{if } \varepsilon_V \ll 1,
    \end{aligned}
    \right.
\end{equation}

By inserting $U$ into the $\theta$ evolution equation \eqref{eqn:the-evo} and using the small-scale expressions \cref{eqn:st-ssd1,eqn:st-ssd2}, we can extract the linear and highest-order terms:
\begin{equation}\label{eqn:st-ssd-full}
    \theta_t = 
    \left\{
    \begin{aligned}
     &\dfrac{\varepsilon_C}{\varepsilon_V}\left(\dfrac{2\pi}{L}\right)^2 \theta_{\alpha\alpha} + \mathscr{N}(\alpha), &\quad \text{if } \varepsilon_V=\mathcal{O}(1), \\
    &2\varepsilon_C\left(\dfrac{2\pi}{L}\right)^3 \mathscr{H}[\theta_{\alpha\alpha\alpha}] + \mathscr{N}(\alpha), &\quad \text{if } \varepsilon_V \ll 1,
    \end{aligned}
    \right.
\end{equation}
where $\mathscr{N}$ consists of the remaining lower-order terms.
The case $\varepsilon_V = \mathcal{O}(1)$ is second-order diffusive and is similar to a heat equation or a mean curvature flow.
The stiffness of the second-order derivative term can be removed by employing an implicit time-stepping scheme for the stiff term.
The second case $\varepsilon_V \ll 1$ is third-order diffusive \cite{Hou1994}.
Implicit schemes for this case are more difficult due to the nonlocal Hilbert transform $\mathscr{H}$.
Using the fact that $\mathscr{H}$ is diagonalizable under the Fourier transform, an accurate and efficient semi-implicit scheme can be devised.

\subsection{Semi-implicit scheme}
The evolution equations \eqref{eqn:L-evo} and \eqref{eqn:ave-evo} are not stiff; they can be discretized with an explicit scheme.
We use the second-order Adams-Bashforth scheme,
\begin{align}
    &L^{n+1} = L^{n} + \tau (3M^n - M^{n-1}), \quad M =  - \int_0^{2\pi} \theta_{\alpha^{\prime}} U \,d\alpha^{\prime},\\
    &\mathbf{\overline{X}}^{n+1} = \mathbf{\overline{X}}^{n} + \tau (3Q^n - Q^{n-1}), \quad Q= \dfrac{1}{2\pi}\int_{0}^{2\pi} U\mathbf{n} + V\mathbf{s}\,d\alpha.
\end{align}
The evolution equation \eqref{eqn:the-evo} involves a stiff term, but it is hidden in the normal velocity $U$ through a DtN mapping.
In order to both remove the stiffness and avoid solving nonlinear algebraic systems, only the highest-order term in \eqref{eqn:the-evo} is discretized implicitly.
We first rewrite equation \eqref{eqn:the-evo} as
\begin{align}\label{eqn:the-evo2}
    \theta_t &=  \dfrac{2\pi}{L} (U_{\alpha}+ \theta_{\alpha}T) + (\mathscr{L}(\alpha) - \mathscr{L}(\alpha))
    = \mathscr{L}(\alpha) + \mathscr{N}(\alpha),
\end{align}
where $\mathscr{L}$ is the highest-order and linear term and $\mathscr{N}$ consists of the remaining lower-order terms.
We set $\mathscr{L}=\mathscr{L}_1=\lambda_1(2\pi/L)^2\theta_{\alpha\alpha}$ for second-order diffusion and  $\mathscr{L}=\mathscr{L}_2=\lambda_2(2\pi/L)^3\mathscr{H}[\theta_{\alpha\alpha\alpha}]$ for third-order diffusion, where $\lambda_1$ and $\lambda_2$ are two constant parameters.
For the Hele-Shaw flow and the Stefan problem with constant $\varepsilon_C$ and $\varepsilon_V$, the parameters $\lambda_1$ and $\lambda_2$ can be chosen to match the corresponding constant factors in the highest-order terms.
For the Stefan problem with anisotropic $\varepsilon_C$ and $\varepsilon_V$, $\lambda_1$ and $\lambda_2$ are regarded as stabilization parameters, which should be chosen large enough to ensure stability.
With a frozen coefficient analysis, we choose the parameters as
\begin{align}
    \lambda_1 = \max_{\alpha\in[0,2\pi)}\left|\dfrac{\varepsilon_C(\mathbf{n}(\alpha))}{\varepsilon_V(\mathbf{n}(\alpha))}\right|, \quad \lambda_2 = \max_{\alpha\in[0,2\pi)} |2\varepsilon_C(\mathbf{n}(\alpha))|.
\end{align}
In the Fourier space, the Hilbert transform $\mathscr{H}$ becomes diagonal, and \eqref{eqn:the-evo2} simplifies to
\begin{align}
    \hat{\theta}_t(k) &= - \lambda_1\left(\dfrac{2\pi}{L}\right)^2 k^2 \hat{\theta}(k) + \mathscr{\hat N}(k)\quad \text{for } \mathscr{L}=\mathscr{L}_1,\label{eqn:Fourier-theta1}\\
    \hat{\theta}_t(k) &= - \lambda_2\left(\dfrac{2\pi}{L}\right)^3 |k|^3 \hat{\theta}(k) + \mathscr{\hat N}(k)\quad \text{for } \mathscr{L}=\mathscr{L}_2.\label{eqn:Fourier-theta2}
\end{align}
A linear propagator and a second-order Adams-Bashforth method are used to discretize the stiff part and the non-stiff part in \eqref{eqn:Fourier-theta1} and \eqref{eqn:Fourier-theta2}, respectively,
\begin{align}
    &\hat{\theta}^{n+1}(k) = e_k(t_n,t_{n+1}) \hat{\theta}^n(k) + \dfrac{\tau}{2}(3e_k(t_n, t_{n+1})\mathscr{\hat N}^n(k) - e_k(t_{n-1}, t_{n+1})\mathscr{\hat N}^{n-1}(k)),
\end{align}
where the factors $e_k(t_n, t_{n+1})$ and $e_k(t_{n-1}, t_{n+1})$ are specified as 
\begin{align}
    &e_k(t_n, t_{n+1})=\exp \left(-\frac{\lambda_1 \tau}{2}(2 \pi k)^2\left[\frac{1}{\left(L^n\right)^2}+\frac{1}{\left(L^{n+1}\right)^2}\right]\right), \\
    &e_k(t_{n-1}, t_{n+1}) = \exp \left(-\lambda_1 \tau(2\pi k)^2 \left[\frac{1}{2\left(L^{n-1}\right)^2}+\frac{1}{\left(L^n\right)^2}+\frac{1}{2\left(L^{n+1}\right)^2}\right]\right),
\end{align}
for the case of $\mathscr{L} = \mathscr{L}_1$ and
\begin{align}
    &e_k(t_n, t_{n+1})=\exp \left(-\frac{\lambda_2 \tau}{2}(2 \pi|k|)^3\left[\frac{1}{\left(L^n\right)^3}+\frac{1}{\left(L^{n+1}\right)^3}\right]\right), \\
    &e_k(t_{n-1}, t_{n+1}) = \exp \left(-\lambda_2 \tau(2\pi|k|)^3 \left[\frac{1}{2\left(L^{n-1}\right)^3}+\frac{1}{\left(L^n\right)^3}+\frac{1}{2\left(L^{n+1}\right)^3}\right]\right),
\end{align}
for the case of $\mathscr{L} = \mathscr{L}_2$.

\begin{remark}
Numerical operations, such as differentiation, integration, and solving ODEs, are all performed in the Fourier space using FFTs. This allows the interface evolution method to have spectral accuracy, resulting in smaller errors compared to second-order time integration methods. It is important to note that $\theta$ is not a smooth periodic function in $\mathbb R/(2\pi\mathbb Z)$. To obtain accurate results, one can use the auxiliary variable $\eta = \theta - \alpha$ when performing Fourier transforms.
\end{remark}
\section{Numerical results}\label{sec:res}
In this section, we assess the proposed method through a series of numerical examples.
All computations were performed on a computer equipped with an Intel(R) Core(TM) i7-10700K CPU at 3.80 GHz and 16 GB of memory, and the numerical experiments were implemented in C++.
We begin with the Hele-Shaw flow, where we study convergence and simulate several bubble dynamics. The method captures complex interface growth and deformation and remains effective in long-time computations that produce intricate finger-like patterns.
We then consider the Stefan problem. Through a range of examples, we examine convergence, stability, and the ability of the method to track solidification interfaces accurately. We also simulate dendritic solidification with and without flow in the liquid phase. Together, these examples illustrate the robustness and versatility of the method across a variety of moving-interface problems.

\subsection{The Hele-Shaw flow}
\subsubsection{Convergence test}
\begin{figure}[htbp]
    \centering
    \subfigure[]{\includegraphics[width=0.45\textwidth]{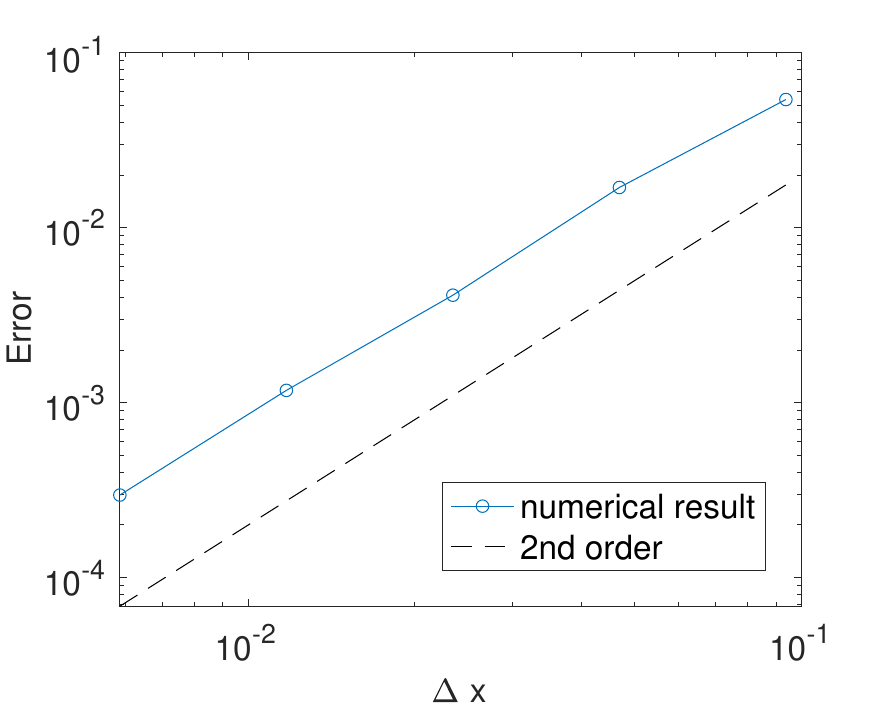}}
    \subfigure[]{\includegraphics[width=0.45\textwidth]{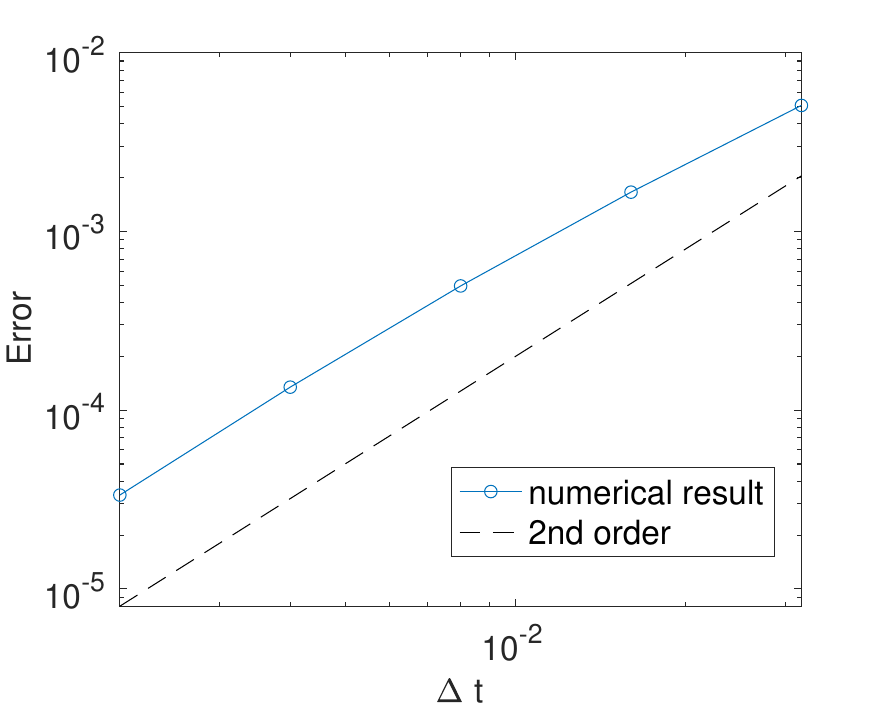}}
    \caption{Spatial and temporal convergence for the Hele-Shaw flow. (a) Spatial accuracy for mesh sizes $\Delta x = 3\times 2^{-l}$, $l = 5,6,\cdots,9$. (b) Temporal accuracy for time steps $\tau = 2^{-l}\times 10^{-3}$, $l = 1,2,\cdots,5$.}
    \label{fig:hs-conv}
\end{figure}

In this example, we investigate the convergence of the method for the Hele-Shaw flow.
The initial shape is a four-fold flower defined by
\begin{align}
    (x(\alpha,0), y(\alpha,0) )= r(\alpha)(\cos\alpha, \sin\alpha), \quad
    r(\alpha) = 0.8 + 0.2\cos 4\alpha, \quad \alpha\in[0, 2\pi).
\end{align}
The surface tension coefficient is $0.01$, and the computation is performed in the bounding box $(-1.5,1.5)^2$.
The air injection rate is set to zero, so the flow preserves area because the fluid is incompressible.
The numerical error is measured by comparing the enclosed area of the curve at $T=1$ with the initial area.
To study the spatial accuracy, we use a fixed time step size of $\tau = 1\times 10^{-3}$ and vary the mesh size as $\Delta x = 3\times 2^{-l}$, with $l = 5,6,\cdots,9$. Similarly, to analyze the temporal accuracy, we use a fixed mesh size of $\Delta x = 3\times 2^{-10}$ and vary the time step size as $\tau = 2^{-l}\times 10^{-3}$, with $l = 1,2,\cdots,5$.
We compute the enclosed area of the curve by numerically evaluating the integral $A[\Gamma] = \frac{1}{2}\int_{\Gamma} \mathbf{x}\cdot\mathbf{n}\,ds$. 
The numerical results are shown in \Cref{fig:hs-conv} and indicate that the method is second-order accurate in both space and time.

\subsubsection{Bubble relaxation}
\begin{figure}[htbp]
    \centering
    \includegraphics[height=0.5\textwidth]{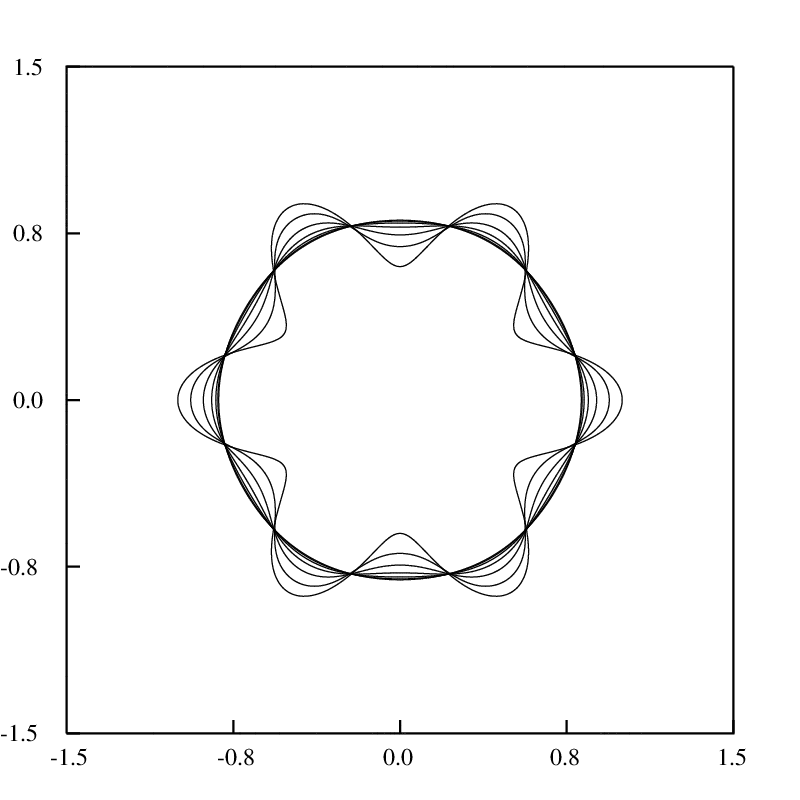}
    \caption{Interface profiles from $t=0$ to $t=1.2$ with time increment $0.2$.}
    \label{fig:hs-stab-1}
\end{figure}

\begin{figure}[htbp]
    \centering
    \subfigure[]{\includegraphics[width=0.45\textwidth]{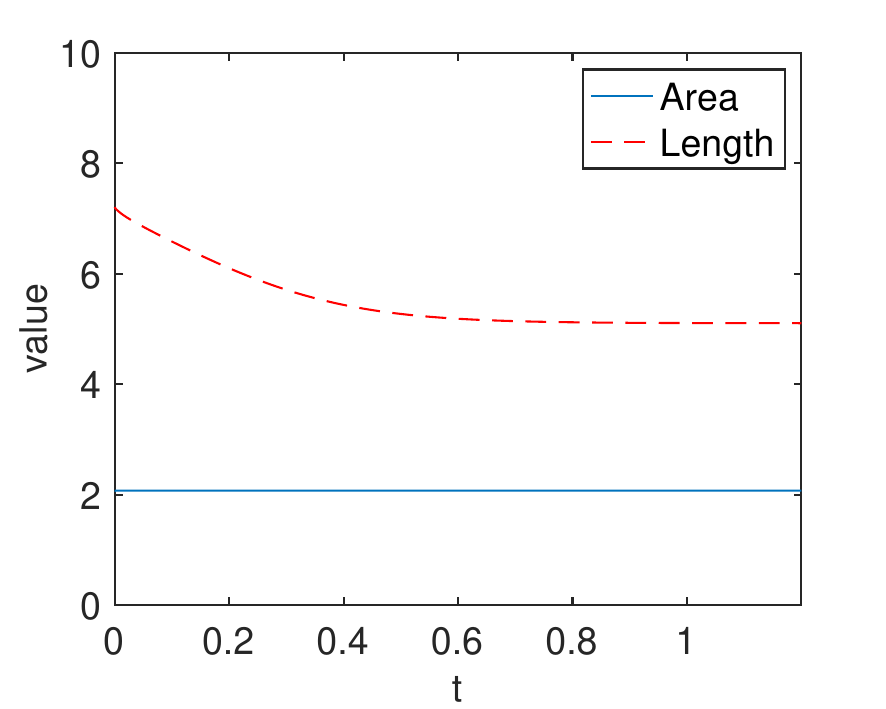}}
    \subfigure[]{\includegraphics[width=0.45\textwidth]{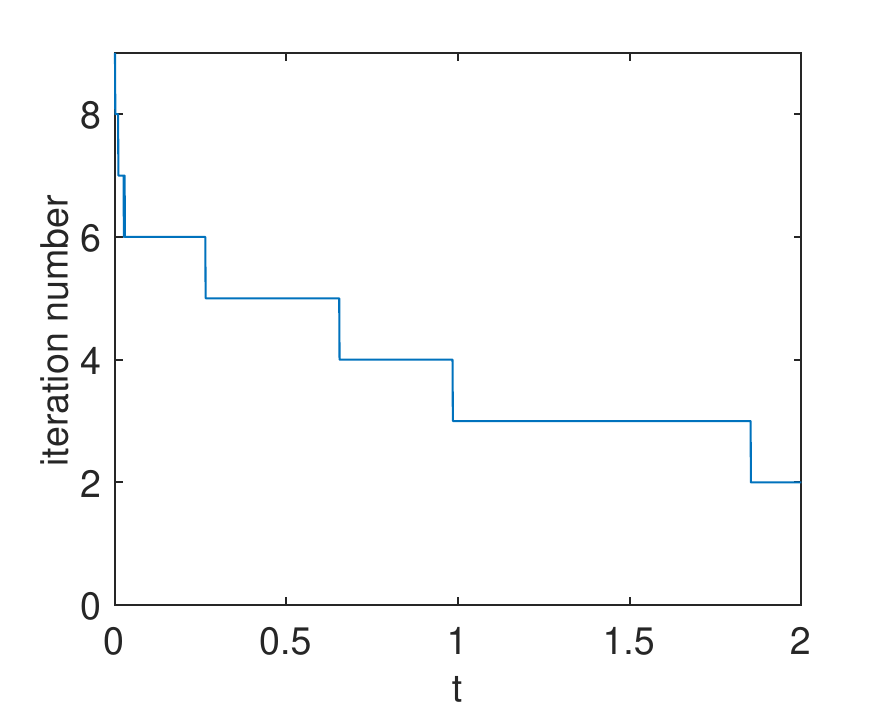}}
    \caption{Bubble relaxation in the Hele-Shaw flow. (a) Time evolution of the enclosed area and interface length. (b) GMRES iteration counts.}
    \label{fig:hs-stab-2}
\end{figure}

The Hele-Shaw flow without air injection exhibits interesting behavior due to the combined effects of surface tension and incompressibility. In particular, the flow is curve-shortening and area-preserving \cite{Sakakibara2021}. To investigate this phenomenon, we consider an initial curve in the form of a six-fold flower, described by the parametric equation
\begin{align}
    (x(\alpha,0), y(\alpha,0) )= r(\alpha)(\cos\alpha, \sin\alpha), \quad
    r(\alpha) = 0.8 + 0.2\cos 6\alpha, \quad \alpha\in[0, 2\pi).
\end{align}
We use a surface tension coefficient of $0.01$ and compute the evolution in the domain $(-1.5,1.5)^2$. The simulation is performed on a $512\times 512$ grid with time step $\tau = 0.001$. The evolution of the curve is shown in \Cref{fig:hs-stab-1}, and the corresponding area and length profiles, together with the GMRES iteration counts, are shown in \Cref{fig:hs-stab-2}.
As a result of the stabilizing influence of surface tension, the initially irregular interface gradually relaxes and approaches a circular shape over time. 
The enclosed area remains constant throughout the evolution, while the interface length decreases and eventually converges.
These observations align with the theoretical understanding of Hele-Shaw flow without air injection.
The GMRES iteration number remains relatively stable and decreases as the curve approaches a circular shape.

\subsubsection{Unstable viscous fingering}
\begin{figure}[htbp]
    \centering
    \subfigure[]{\includegraphics[width=0.4\textwidth]{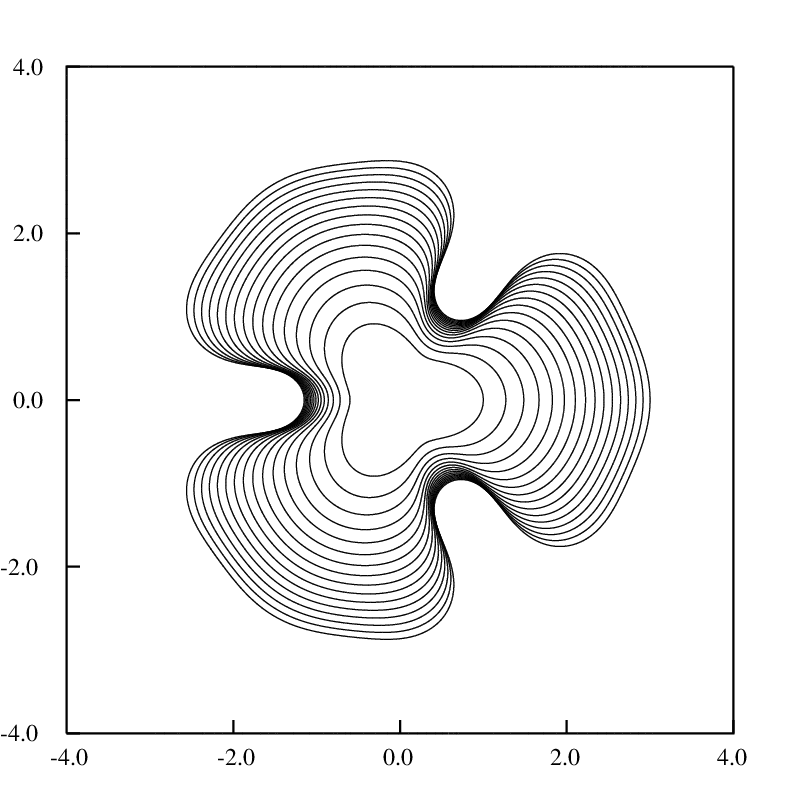}}
    \subfigure[]{\includegraphics[width=0.4\textwidth]{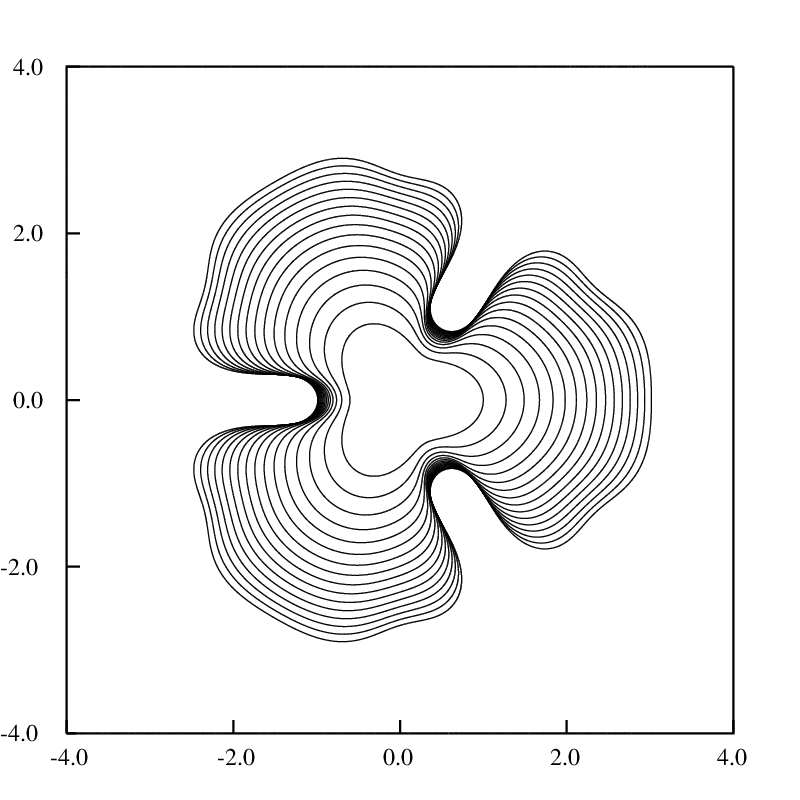}}
    \subfigure[]{\includegraphics[width=0.4\textwidth]{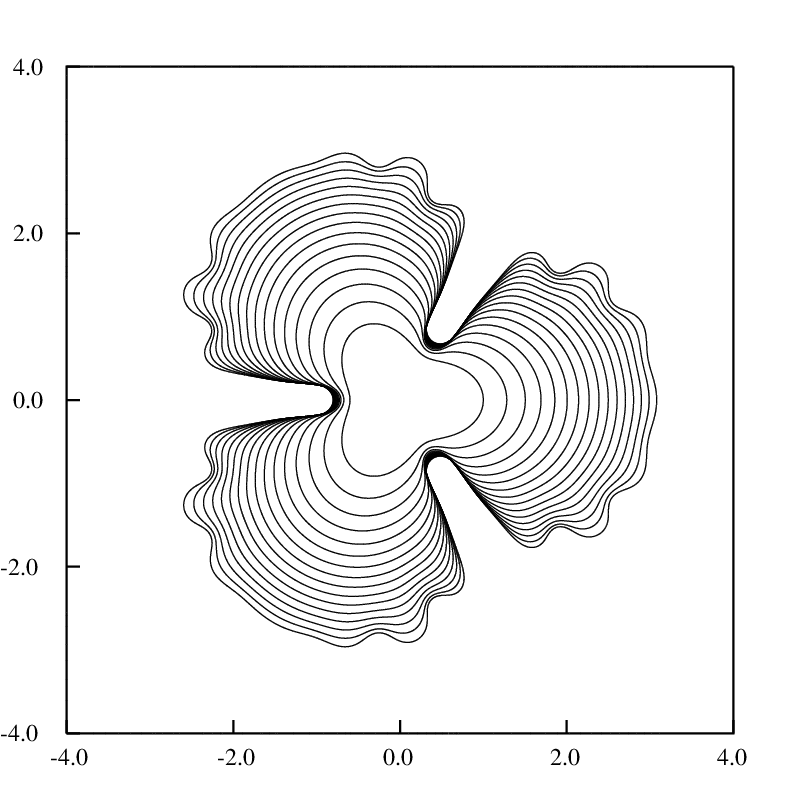}}
    \subfigure[]{\includegraphics[width=0.4\textwidth]{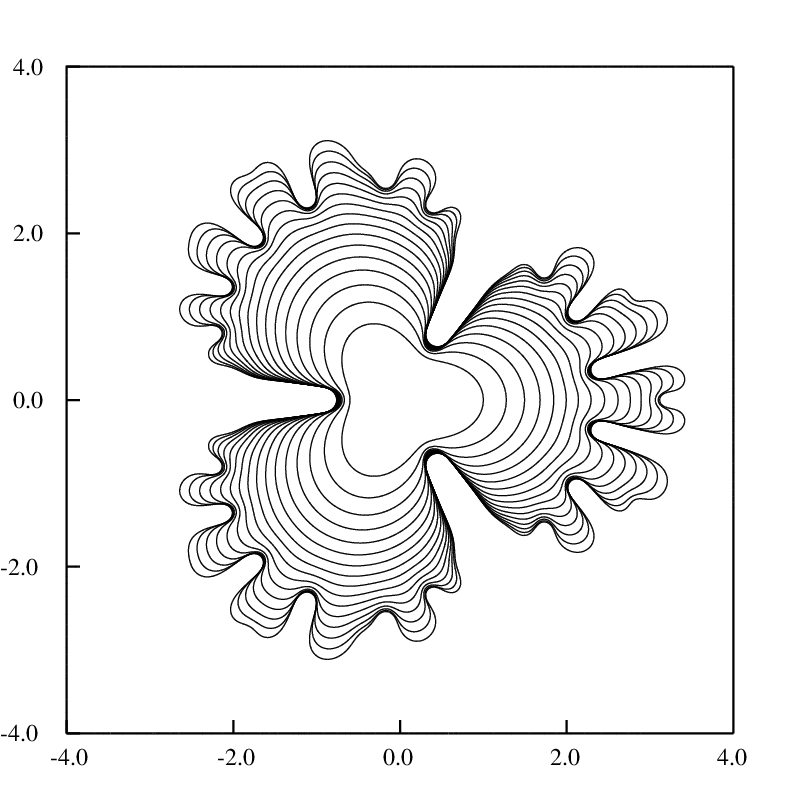}}
    \caption{Interface profiles for unstable viscous fingering with different surface tension coefficients: (a) $1\times 10^{-2}$, (b) $5\times 10^{-3}$, (c) $1\times 10^{-3}$, and (d) $5\times 10^{-4}$.}
    \label{fig:hs-sur}
\end{figure}

The initial curve is given by a three-fold flower,
\begin{align} 
    (x(\alpha,0), y(\alpha,0) )= r(\alpha)(\cos\alpha, \sin\alpha), \quad
    r(\alpha) = 0.8 + 0.2\cos 3\alpha, \quad \alpha\in[0, 2\pi).
\end{align}
The air injection rate is set to $J = 1$, so the air bubble grows unstably.
The surface tension coefficient varies from $1\times 10^{-2}$ to $5\times 10^{-4}$.
The computational domain is $(-4,4)^2$.
We solve the problem on a $512\times 512$ grid with time step $\tau = 0.005$.
The interface profiles are shown in \Cref{fig:hs-sur} at time intervals of $0.2$ up to the final time $T=3$.
The competition between the stabilizing effect of surface tension and the destabilizing effect of the driving force due to the injection leads to the viscous fingering feature of the Hele-Shaw flow.
With small surface tension, which has a less stabilizing effect, the growing bubble develops more branches.
For small surface tensions, the symmetry of the interface becomes harder to preserve because grid-induced anisotropy has a stronger influence on the evolution.

\subsubsection{Long-time computation}
\begin{figure}[htbp]
    \centering
    \subfigure[]{\includegraphics[width=0.5\textwidth]{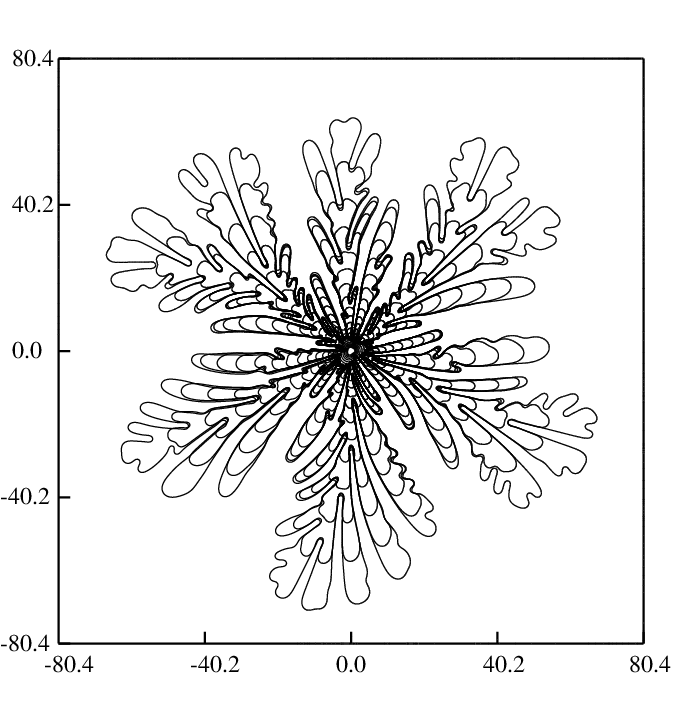}}\\
    \subfigure[]{\includegraphics[width=0.32\textwidth]{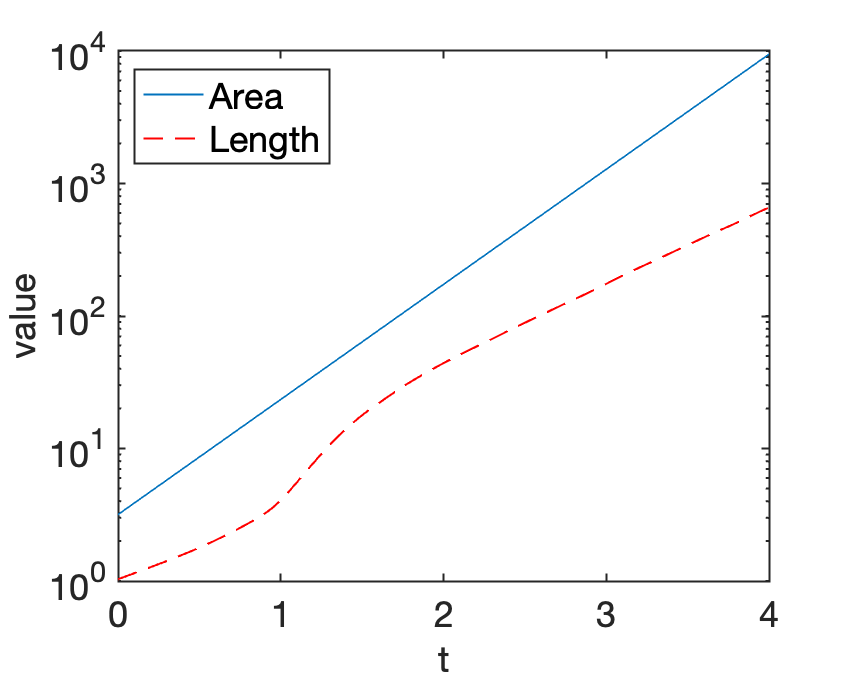}}
    \subfigure[]{\includegraphics[width=0.32\textwidth]{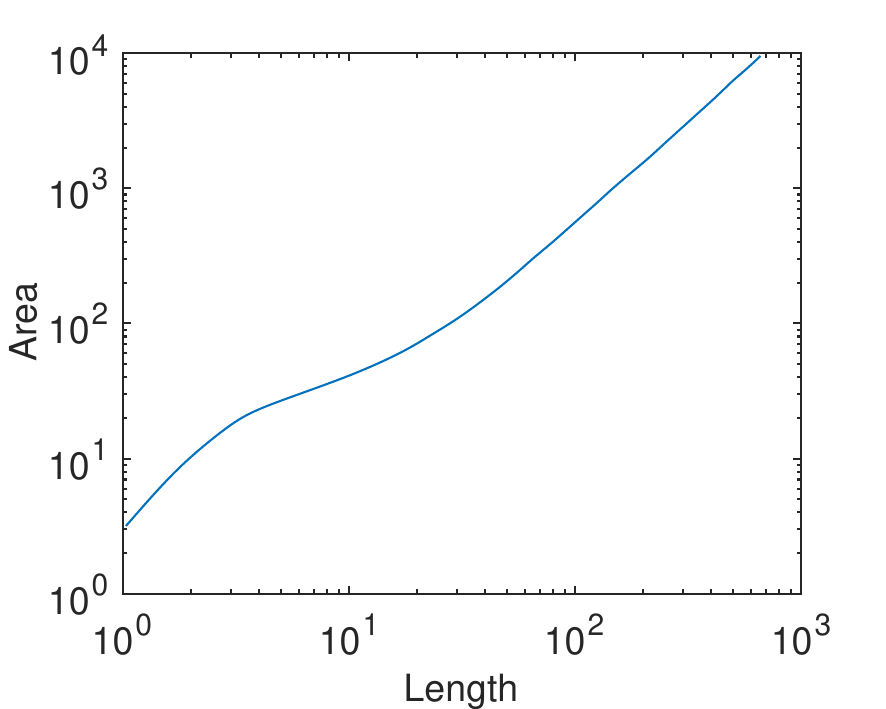}}
    \subfigure[]{\includegraphics[width=0.32\textwidth]{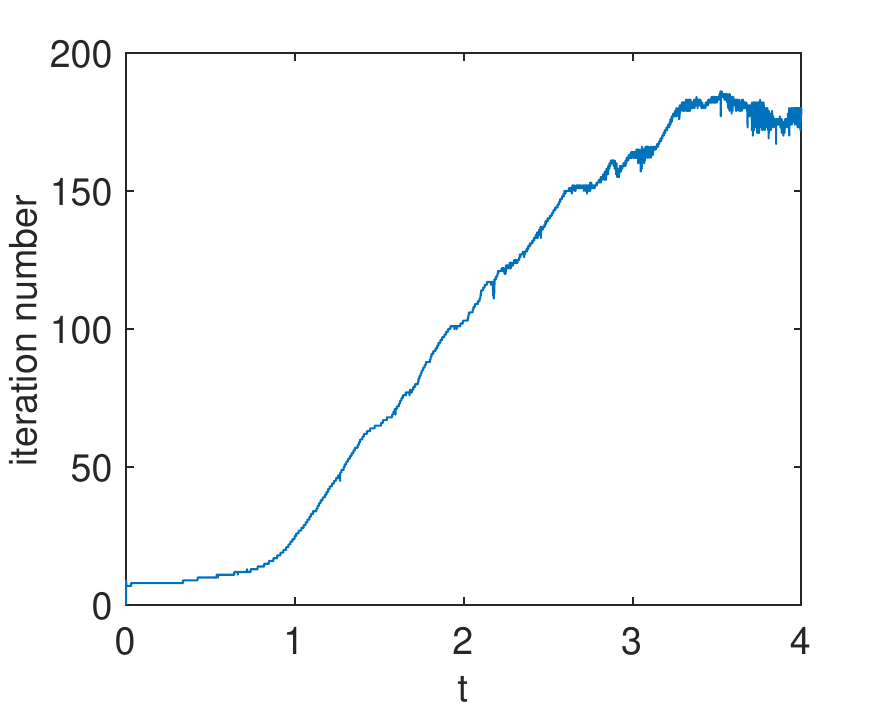}}
    \caption{Long-time Hele-Shaw simulation. (a) Interface histories. (b) and (c) Time evolution of the interface length and enclosed area. (d) GMRES iteration counts.}
    \label{fig:hs-large}
\end{figure}

In this example, we perform a long-time simulation using the present method together with the spatiotemporal rescaling scheme in \cite{Li2007} to compute a large Hele-Shaw bubble.
The computation is carried out in the scaled frame and then mapped back to the original frame.
The initial shape is taken to be a nucleus:
\begin{align}
    (\bar x(\alpha,0), \bar y(\alpha,0) )= \bar r(\alpha)(\cos\alpha, \sin\alpha), \quad
    \bar r(\alpha) = 1.0 + 0.1(\sin 2\alpha + \cos 3\alpha), \quad \alpha\in[0, 2\pi).
\end{align}
The computational domain is $(-1.7,1.7)^2$.
We use a $1024\times 1024$ grid and a time step of $\Delta \bar t = 2\times 10^{-4}$.
The interface points are adaptively refined with the criteria $\Delta \bar s > 1.5 \Delta \bar x$.
The injection rate and the surface tension coefficient are set as $\bar J=1$ and $\sigma=0.001$, respectively.
The interface profiles in the original frame are shown in \Cref{fig:hs-large} at scaled time intervals of $0.2$.
The computation takes $5$ hours to reach the final scaled time $\bar T=3$ (corresponding to the unscaled time $T=203$).
At the final time, the number of marker points on the interface is refined to $16384$, and the enclosed area and the length of the interface are $A=\bar R^2 \bar A=1280$ and $L = \bar R \bar L/2\pi = 169$.

\subsection{The Stefan problem}
\subsubsection{Grid refinement analysis}
\begin{figure}[htbp]
    \centering
    \subfigure[]{\includegraphics[width=0.4\textwidth]{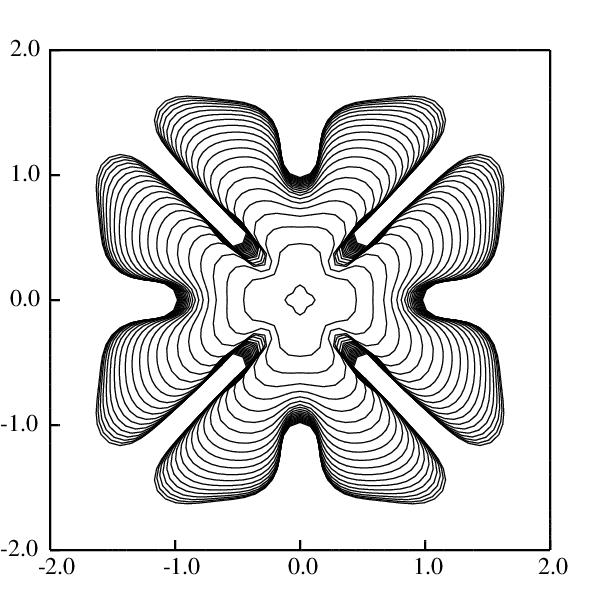}}
    \subfigure[]{\includegraphics[width=0.4\textwidth]{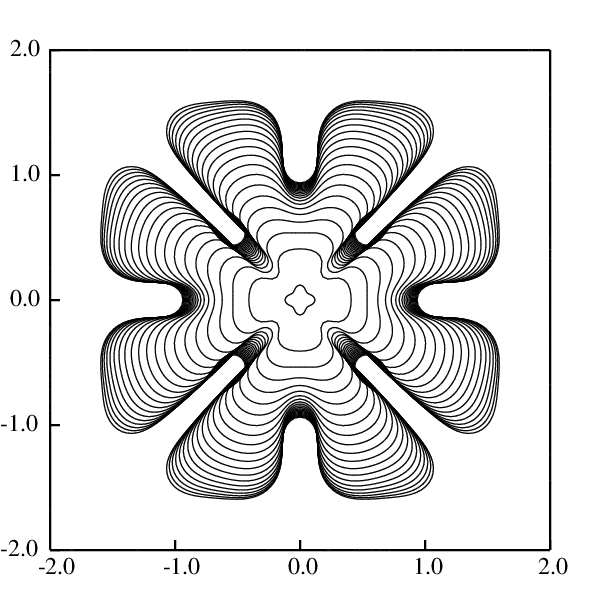}}
    \subfigure[]{\includegraphics[width=0.4\textwidth]{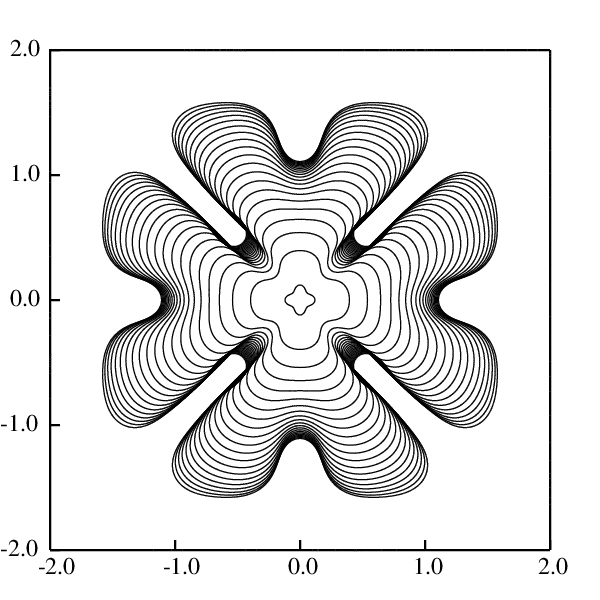}}
    \subfigure[]{\includegraphics[width=0.4\textwidth]{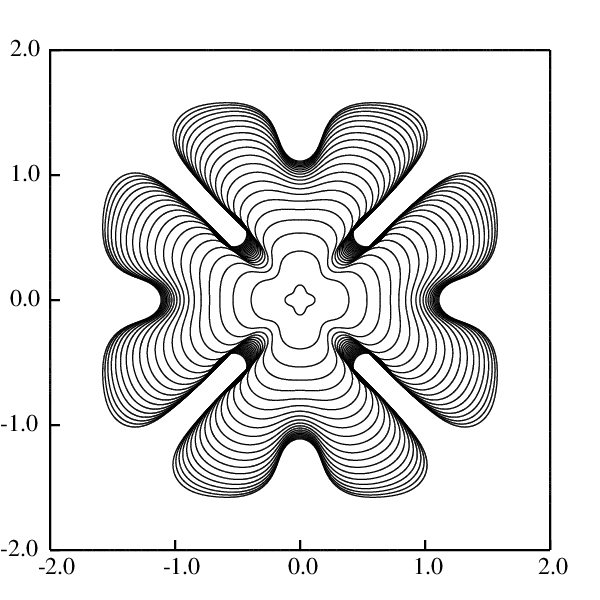}}
    \caption{Grid-refinement study for the Stefan problem: (a) $64\times64$, (b) $128\times128$, (c) $256\times256$, and (d) $512\times512$.}
    \label{fig:st-grid}
\end{figure}

In this example, we present a benchmark grid-refinement study for the Stefan problem.
Initially, a solid seed is placed in an undercooled liquid domain $\mathcal{B}=(-2,2)^2$.
The initial shape is a four-fold flower, which is given by
\begin{equation}
    (x(\alpha,0), y(\alpha,0) )= r(\alpha)(\cos\alpha, \sin\alpha), \quad
    r(\alpha) = 0.1 + 0.02\cos 4\alpha, \quad \alpha\in[0, 2\pi).
\end{equation}
The initial undercooling is $St = -0.5$.
The isotropic surface tension and kinematic coefficients are set to $\varepsilon_C = \varepsilon_V = 2\times 10^{-3}$.
We take the time step as $\tau = 0.001$ and successively refine the grid from $64$ to $512$.
The liquid-solid interface profiles at time intervals of $0.05$ up to the final time $T=0.8$ are shown in \Cref{fig:st-grid}.
Compared with previous results obtained by the level-set and front-tracking methods \cite{Juric1996,Chen1997}, our method exhibits less grid-induced anisotropy and therefore better convergence.
The method is also better able to preserve the symmetries of the interface, even with a coarse grid.

\subsubsection{Stability test}
\begin{figure}[htbp]
    \centering
    \subfigure[]{\includegraphics[width=0.4\textwidth]{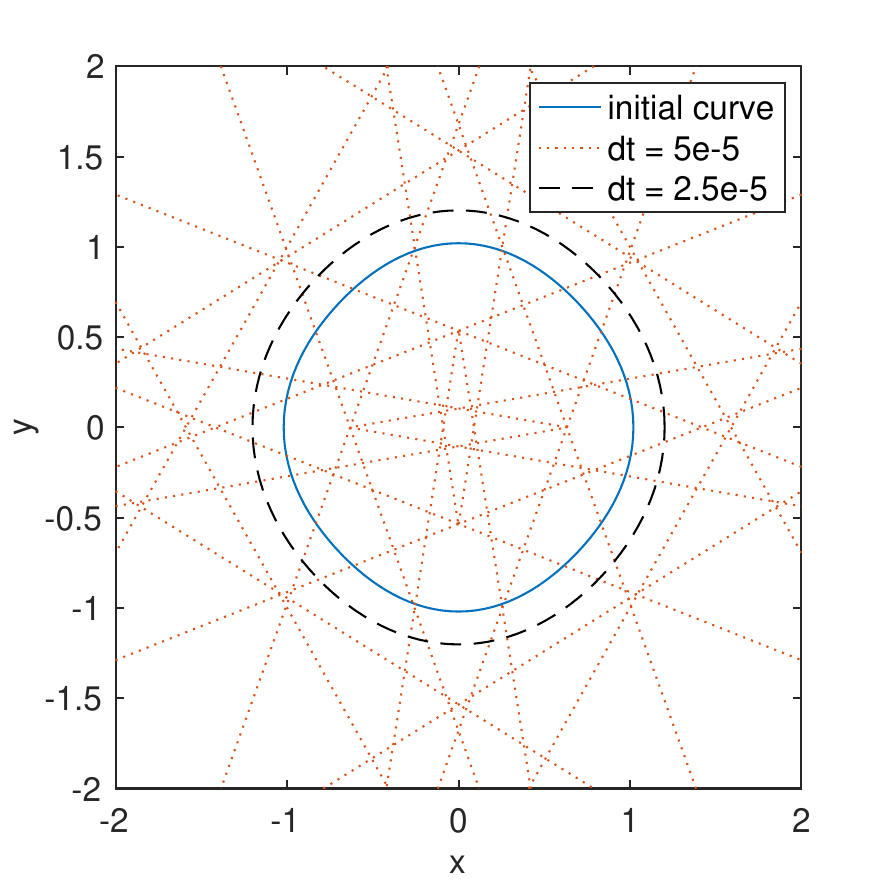}\label{fig:stefan-unstable}}
    \subfigure[]{\includegraphics[width=0.4\textwidth]{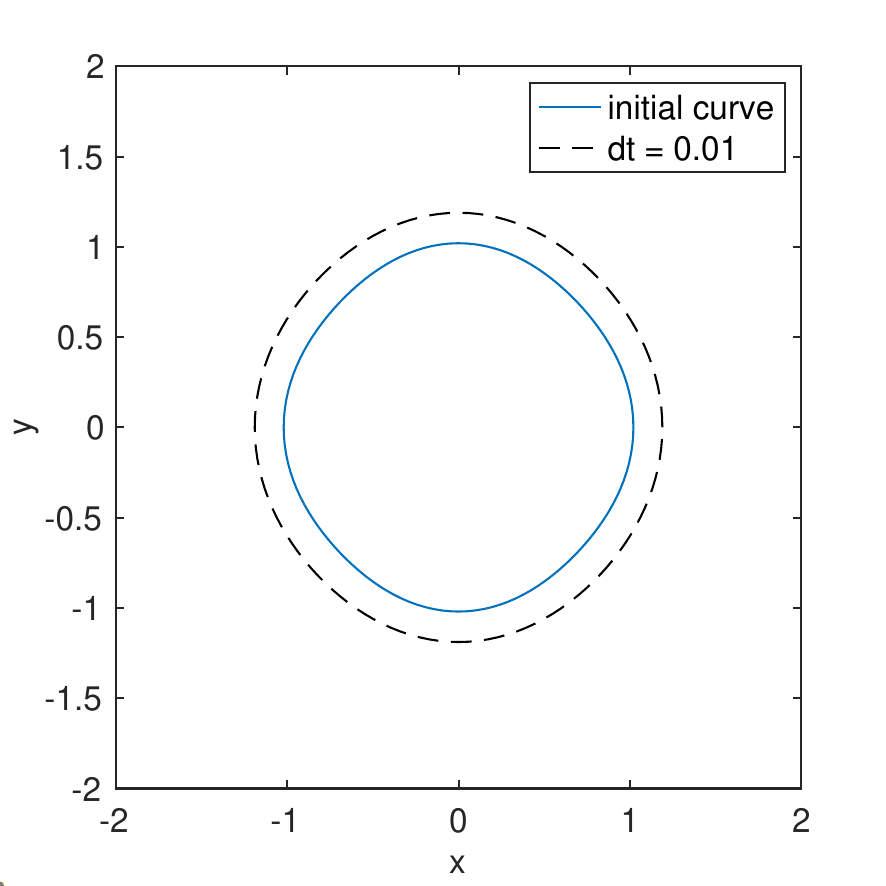}\label{fig:stefan-stable}}
    \caption{Stability test for the Stefan problem with different time-stepping methods: (a) Adams-Bashforth and (b) semi-implicit.}
    \label{fig:st-stab}
\end{figure}

To demonstrate stability, we compare the semi-implicit scheme with an explicit Adams-Bashforth scheme.
We consider the case $\varepsilon_V = 0$ and $\varepsilon_C = 0.05$, which leads to third-order stiffness.
The initial shape is a slightly perturbed circle
\begin{equation}
    (x(\alpha,0), y(\alpha,0) )= r(\alpha)(\cos\alpha, \sin\alpha), \quad
    r(\alpha) = 1 + 0.02\cos 4\alpha, \quad \alpha\in[0, 2\pi).
\end{equation}
We use a $256\times 256$ grid for the computational domain $\mathcal{B}=(-2,2)^2$. 
The interface is discretized with $128$ points.
First, we solve the problem with the explicit scheme using time steps $\tau = 5\times10^{-5}$ and $\tau = 2.5\times 10^{-5}$.
The interface profiles at $t=0.1$ are shown in \Cref{fig:stefan-unstable}.
The smaller time step $\tau = 2.5\times 10^{-5}$ is stable, whereas the larger one $\tau = 5\times 10^{-5}$ is unstable and the solution quickly blows up.
We then solve the same problem with the semi-implicit scheme using time step $\tau = 0.01$, and the solution at $t=0.1$ is shown in \Cref{fig:stefan-stable}.
Compared with the explicit scheme, the semi-implicit scheme remains stable with a much larger time step and is therefore far more efficient.

\subsubsection{Comparison with the solvability theory}
\begin{figure}[htbp]
    \centering
    \subfigure[]{\includegraphics[height=0.4\textwidth]{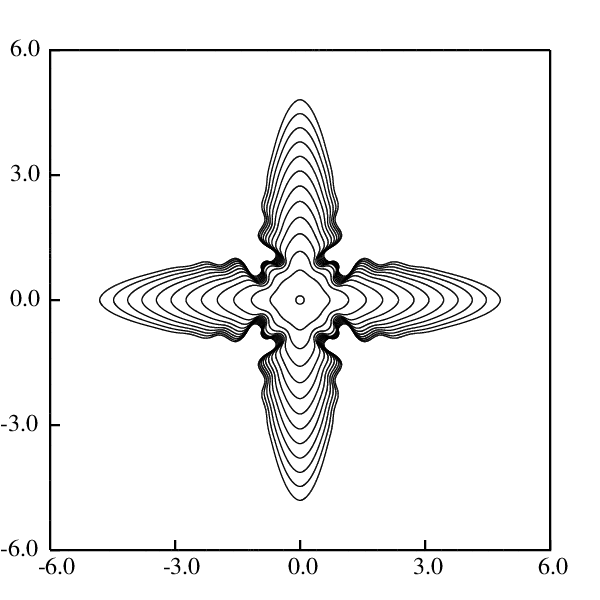}}
    \subfigure[]{\includegraphics[height=0.4\textwidth]{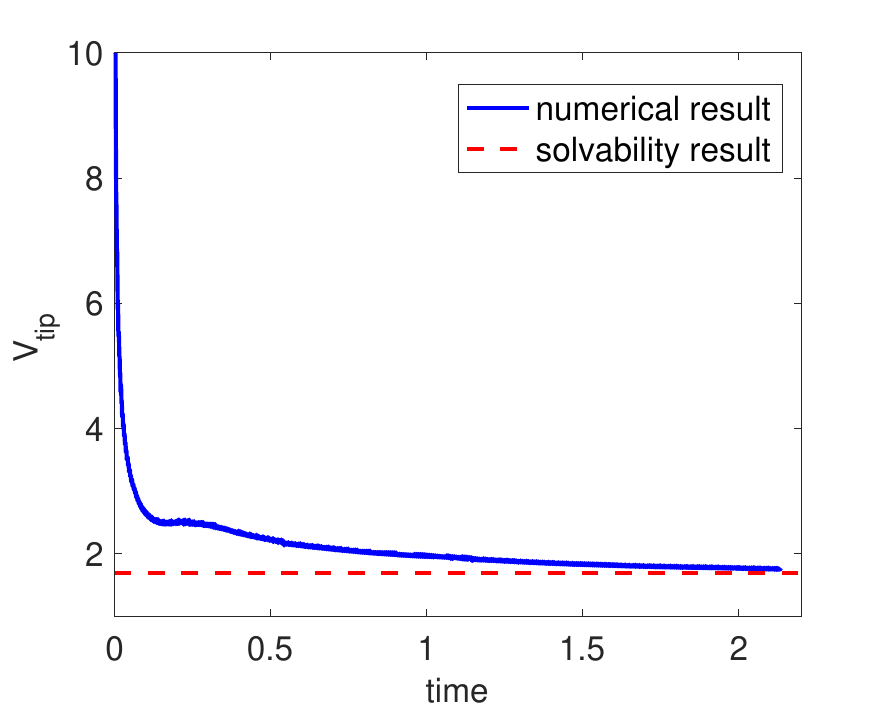}}
    \caption{Comparison with solvability theory for dendritic growth: (a) interface histories and (b) time evolution of the tip velocity.}
    \label{fig:st-solv}
\end{figure}

In a dendritic growth problem, the dendrite growth rate can be predicted by solvability theory \cite{Meiron1986,Chen2009}.
We compare the numerical result with the theoretical prediction to assess the accuracy of the method.
In this example, the initial seed is a circle of radius $0.1$.
The parameters in the Gibbs-Thomson relation are chosen as $\varepsilon_V = 0$ and $\varepsilon_C(\alpha) = 0.001[1+0.4(1-\cos(4\alpha))]$, where $\alpha$ is the angle between the interface normal and the $x$-axis.
The computational domain is $(-6,6)^2$.
We set $\tau = 0.001$ for the computation.
The liquid-solid interface profiles and the tip velocity from $t=0$ to $t=2.2$ are shown in \Cref{fig:st-solv}.
The tip velocity converges to a value consistent with the solvability prediction of $1.7$.

\subsubsection{Anisotropic dendritic growth}
\begin{figure}[htbp]
    \centering
    \subfigure[Numerical results with $\alpha_0 = 0$.]{
    \includegraphics[width=0.24\textwidth]{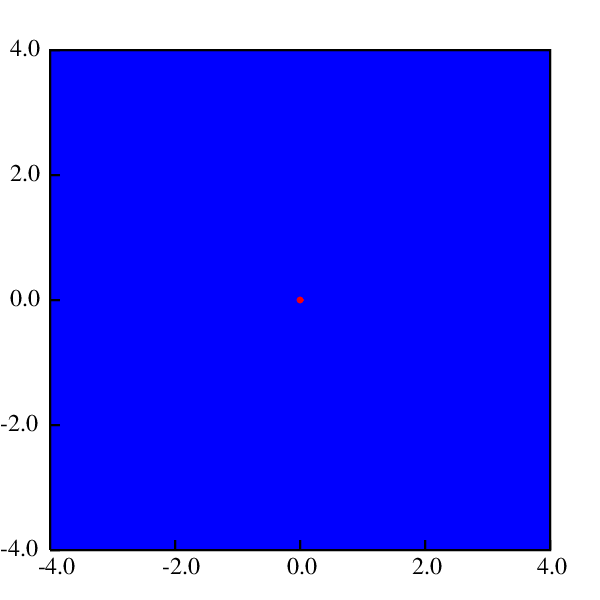}
    \includegraphics[width=0.24\textwidth]{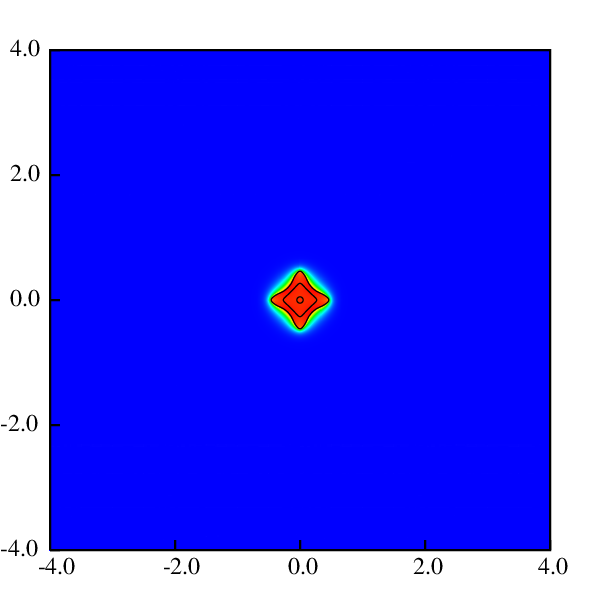}
    \includegraphics[width=0.24\textwidth]{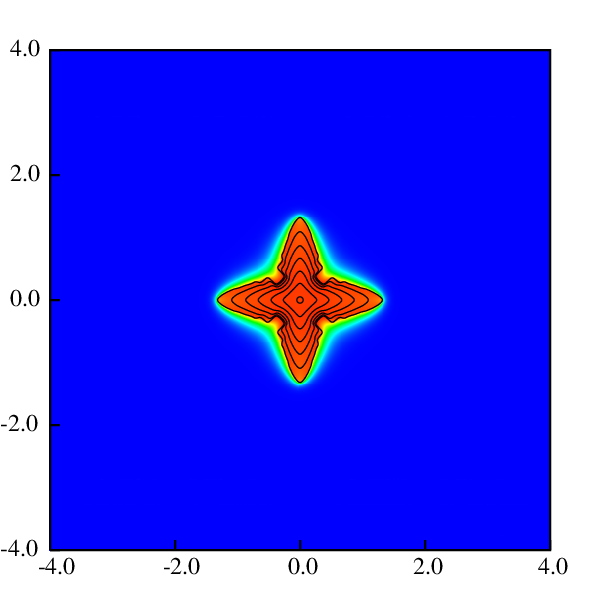}
    \includegraphics[width=0.24\textwidth]{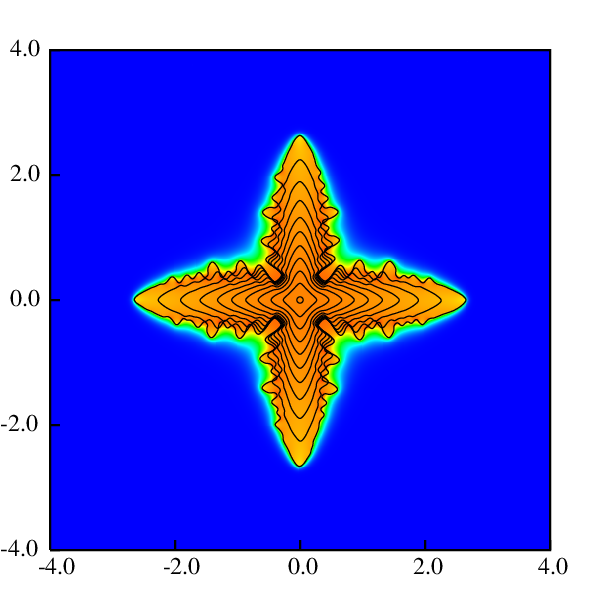}}
    \subfigure[Numerical results with $\alpha_0=\pi/4$.]{
    \includegraphics[width=0.24\textwidth]{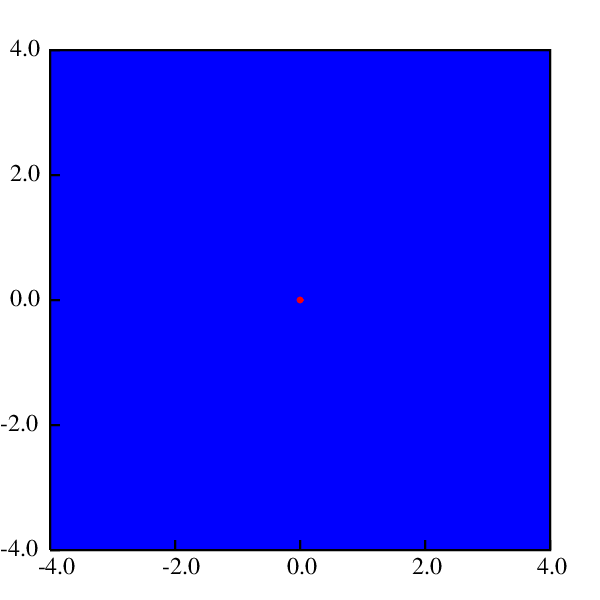}
    \includegraphics[width=0.24\textwidth]{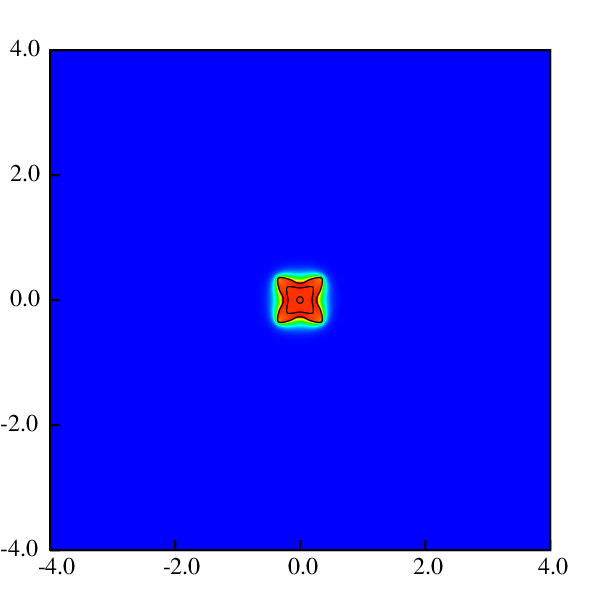}
    \includegraphics[width=0.24\textwidth]{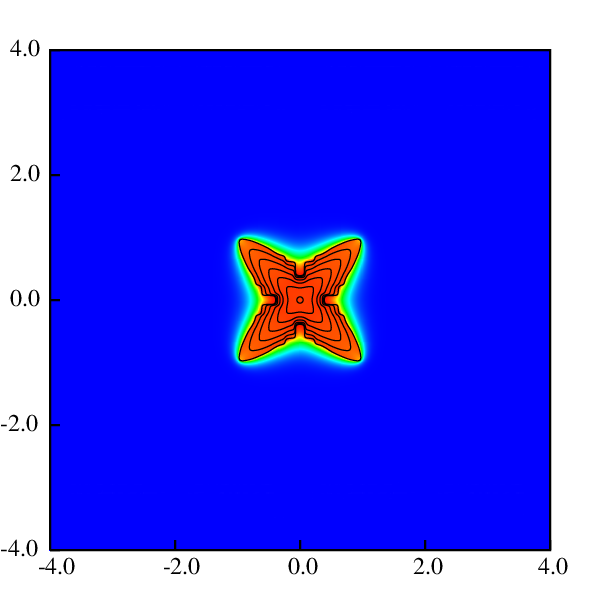}
    \includegraphics[width=0.24\textwidth]{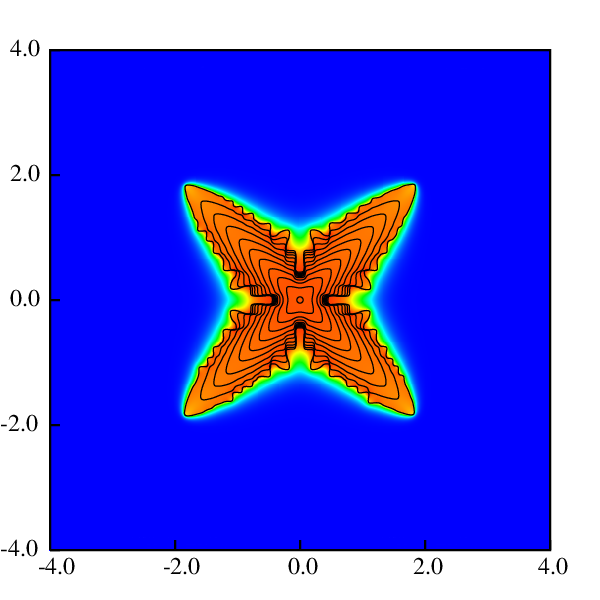}}
    \caption{Dendritic growth with four-fold anisotropy. Snapshots of the interface and temperature field are shown at $t=0$, $0.02$, $0.06$, and $0.1$.}
    \label{fig:st-aniso-4}
\end{figure}

We consider two anisotropy cases. In the first, we use a four-fold anisotropy with $\varepsilon_V = 0.002$ and $\varepsilon_C(\alpha) = 0.002(8/3\sin^4(2(\alpha-\alpha_0)))$, where $\alpha_0=0$ or $\pi/4$. The initial seed is a circle of radius $0.05$ centered in the computational domain $(-4,4)^2$, which is filled with undercooled liquid at $St = -0.65$.
The interface and temperature field are shown in \Cref{fig:st-aniso-4}.

\begin{figure}[htbp]
    \centering
    \subfigure[Numerical results with $St=-0.55$.]{
    \includegraphics[width=0.24\textwidth]{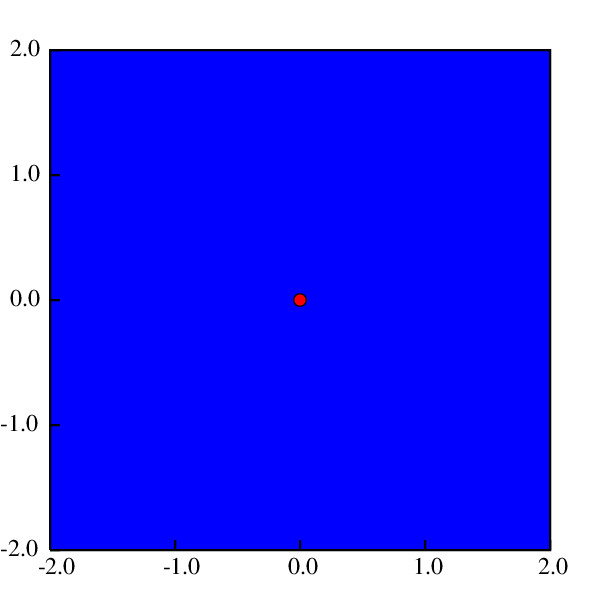} 
    \includegraphics[width=0.24\textwidth]{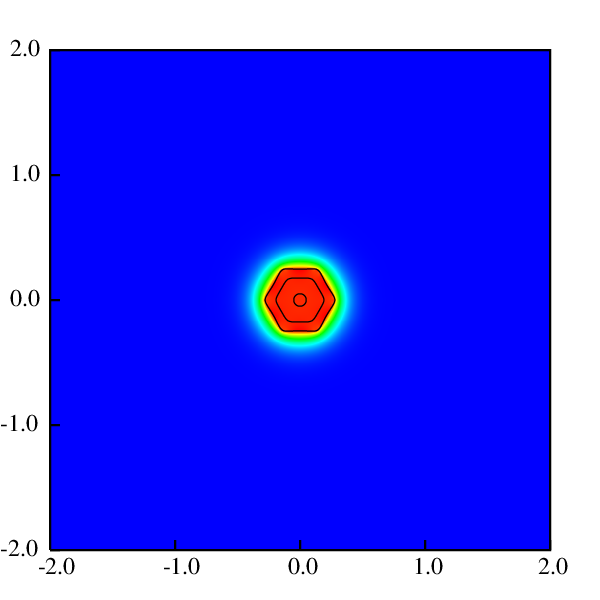} 
    \includegraphics[width=0.24\textwidth]{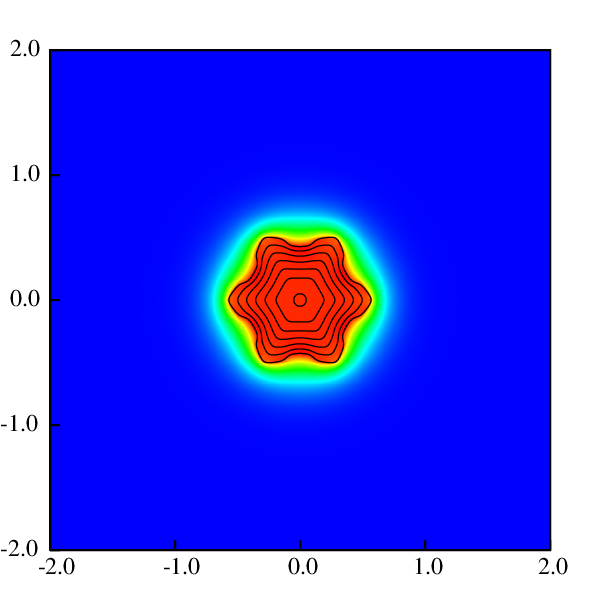} 
    \includegraphics[width=0.24\textwidth]{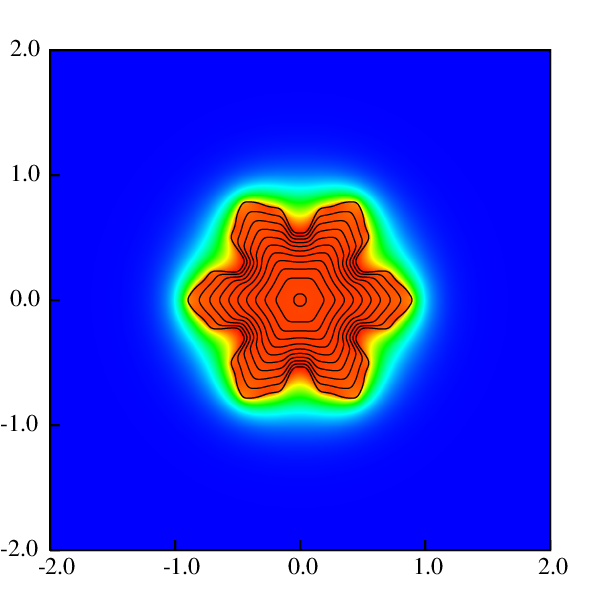}}
    \subfigure[Numerical results with $St=-0.65$.]{
    \includegraphics[width=0.24\textwidth]{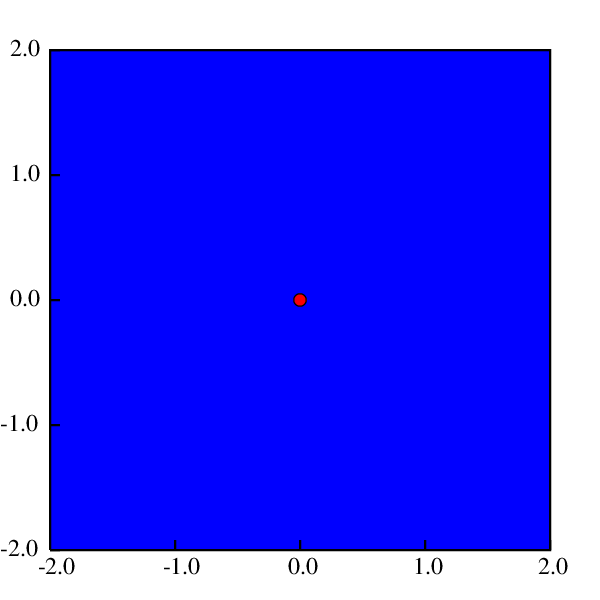} 
    \includegraphics[width=0.24\textwidth]{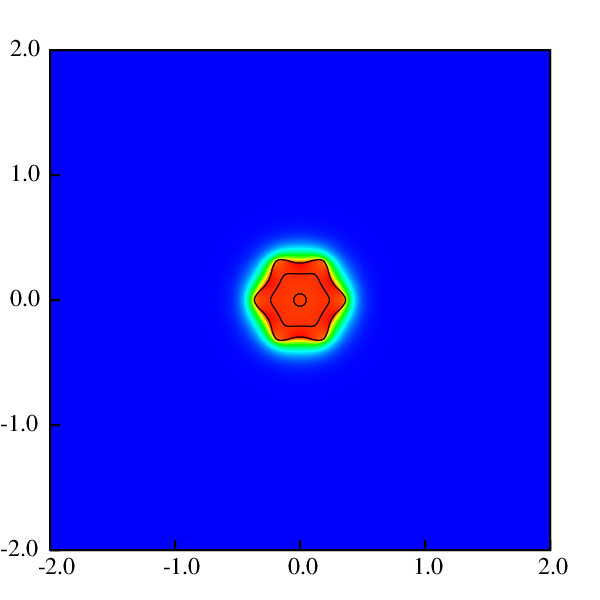} 
    \includegraphics[width=0.24\textwidth]{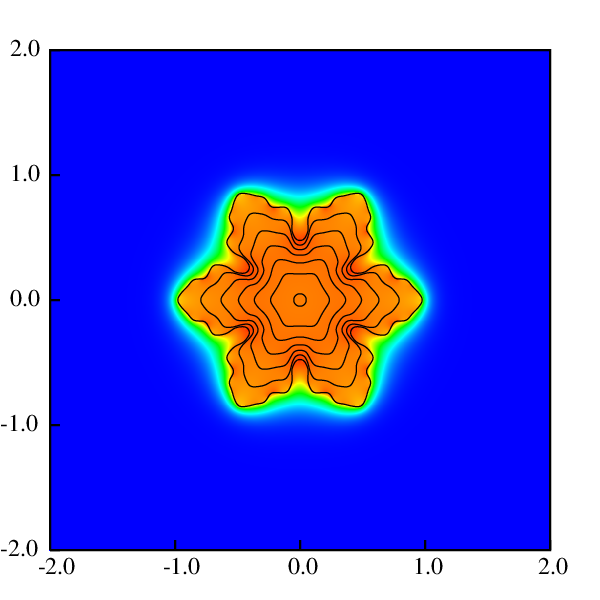} 
    \includegraphics[width=0.24\textwidth]{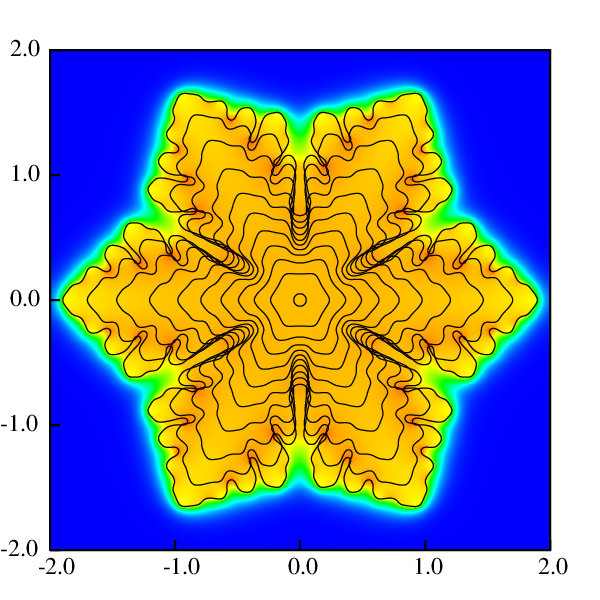}}
    \caption{Interface morphologies and the temperature field of the dendritic growth problem with six-fold anisotropy. Snapshots are taken at $t=0, 0.02, 0.06$, and $0.1$.}
    \label{fig:st-snowflake}
\end{figure}

In the second case, we consider a six-fold anisotropy with $\varepsilon_V = 0.002$ and $\varepsilon_C(\alpha) = 0.002(8/3\sin^4(3\alpha))$. The initial seed is again a circle of radius $0.05$ centered in the computational domain $\mathcal{B} = (-2,2)^2$. We consider two undercooling numbers, $St=-0.55$ and $St=-0.65$. The resulting six-fold anisotropy produces snowflake-shaped dendritic patterns, shown in \Cref{fig:st-snowflake}.

\subsubsection{Dendritic growth with external flow} \label{eg:conv}
\begin{figure}[htbp]
    \centering
    \includegraphics[width=0.24\textwidth]{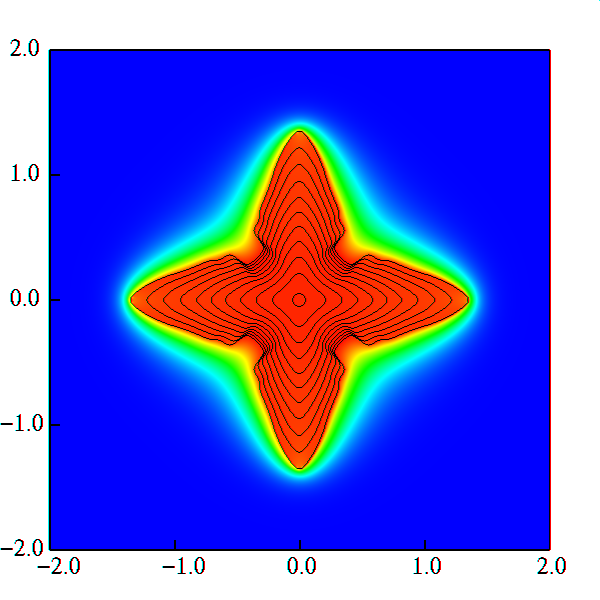} 
    \includegraphics[width=0.24\textwidth]{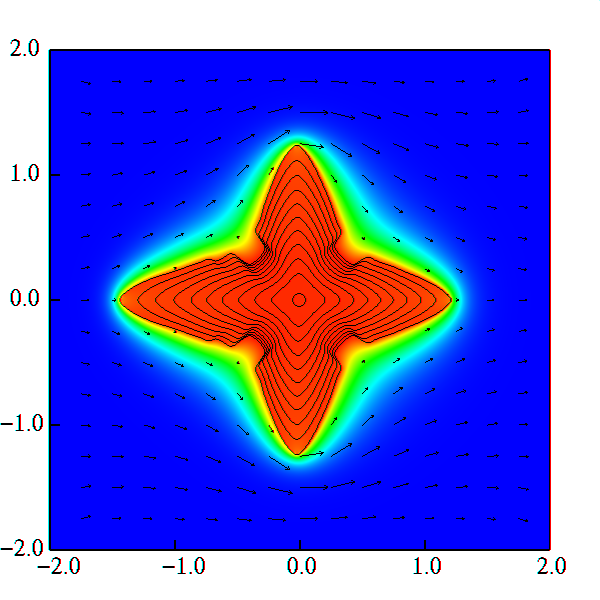} 
    \includegraphics[width=0.24\textwidth]{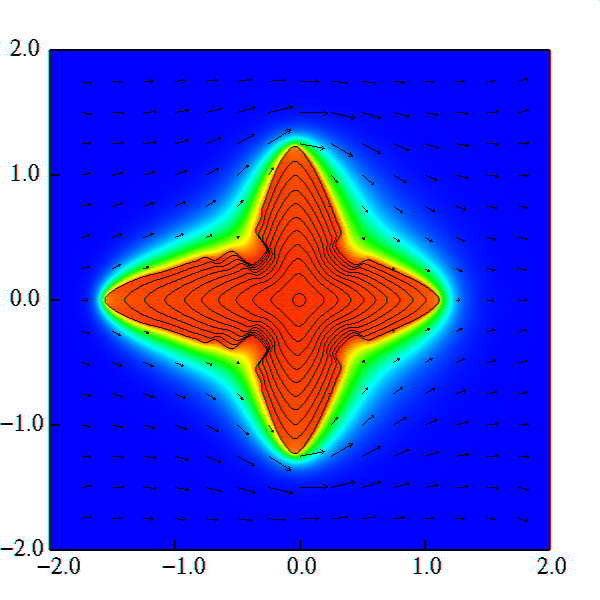} 
    \includegraphics[width=0.24\textwidth]{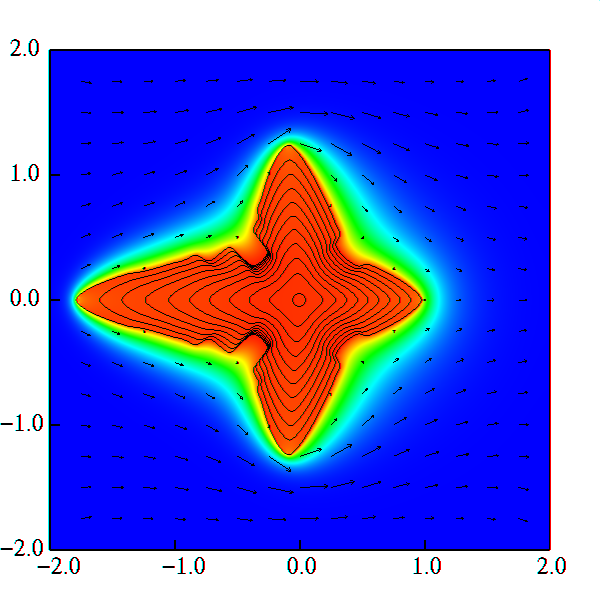} 
    \caption{Dendritic growth with imposed flow. Interface shape, temperature field, and flow field at $t=0.1$ for inflow velocities $u_0 = 0$, $2$, $4$, and $8$ from left to right.}
    \label{fig:convection}
\end{figure}

\begin{figure}[htbp]
    \centering
    \includegraphics[width=0.4\textwidth]{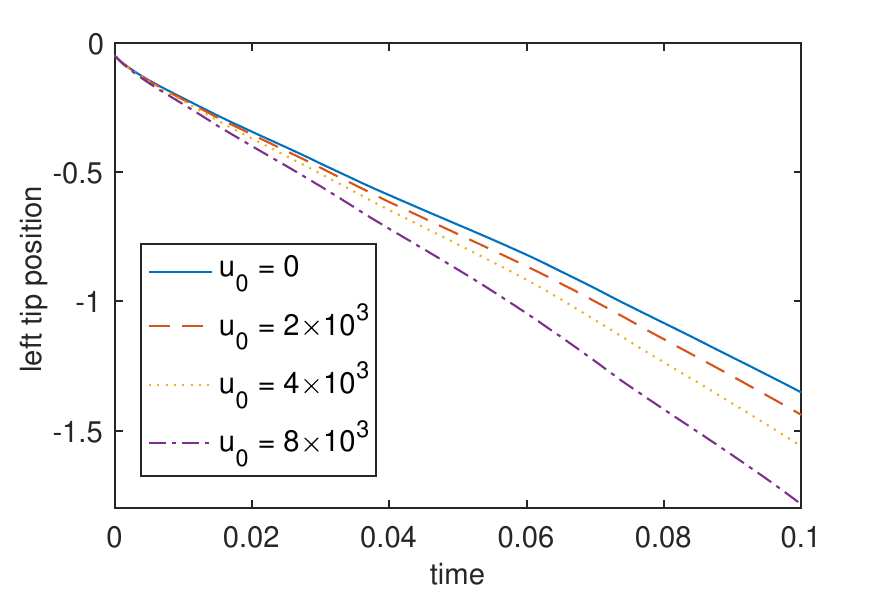} 
    \includegraphics[width=0.4\textwidth]{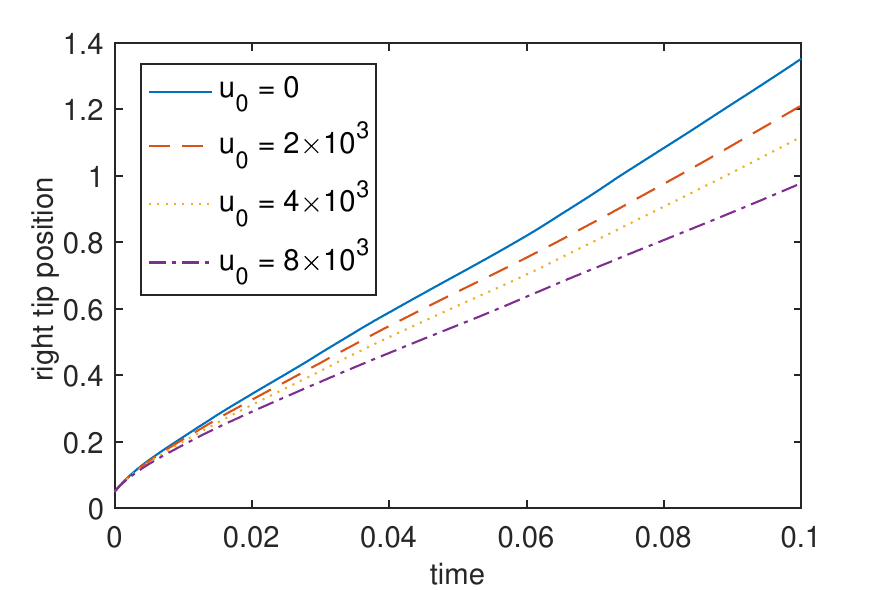} 
    \caption{Time evolution of the $x$-coordinates of the left and right tips for different flow velocities.}
    \label{fig:conv-tip}
\end{figure}

We next examine convection effects in dendritic growth.
In this example, we solve the Stefan problem with flow in the liquid phase.
Initially, the seed is a circle of radius $0.05$ centered in the computational domain $(-2,2)^2$ and surrounded by undercooled fluid with temperature $St = -0.5$.
We use four-fold anisotropic surface tension $\varepsilon_C(\alpha) = 0.002(8/3\sin^4(2\alpha))$ and local kinematic equilibrium $\varepsilon_V = 0$.
The computation is performed with a $512\times 512$ grid and a time step $\tau = 0.0002$.
Inflow and outflow boundary conditions, $\mathbf{u} = (u_0,0)^T$, are imposed on the left and right boundaries, respectively.
No-slip boundary conditions, $\mathbf{u}=\mathbf{0}$, are imposed on the top and bottom boundaries.
In this example, we neglect buoyancy, that is, $\beta = 0$.
As a result, the flow is driven only by the boundary conditions.
When $u_0 = 0$, the problem is identical to the classical Stefan problem without natural convection.
To examine the effect of convection on the growth pattern, we compare solutions with different flow velocities in the liquid phase.
\Cref{fig:convection} shows the dendritic shape together with the temperature and flow fields at $t=0.1$.
The evolution histories of the $x$-components of the left and right tips are shown in \Cref{fig:conv-tip}.
Convection leads to faster growth of the left branch and slower growth of the right branch.
This effect is more evident as the flow velocity increases.
The released latent heat is transported from left to right, which creates an asymmetric temperature distribution in the $x$-direction.


\subsubsection{Dendritic growth with buoyancy-driven flow}
\begin{figure}[htbp]
    \centering
    \includegraphics[width=0.24\textwidth]{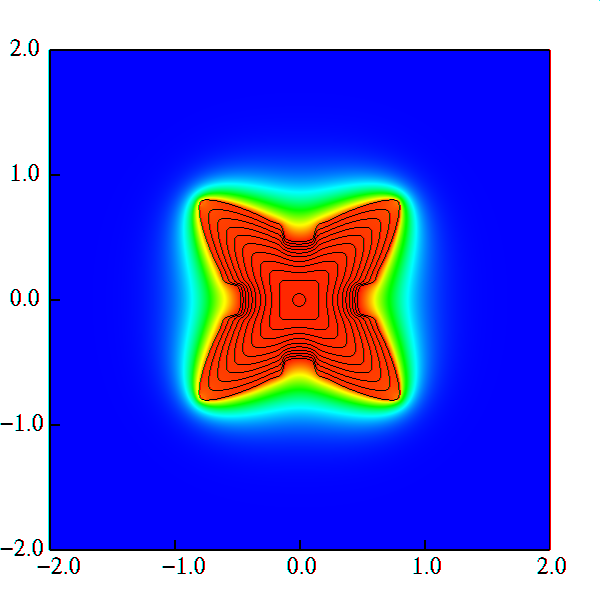} 
    \includegraphics[width=0.24\textwidth]{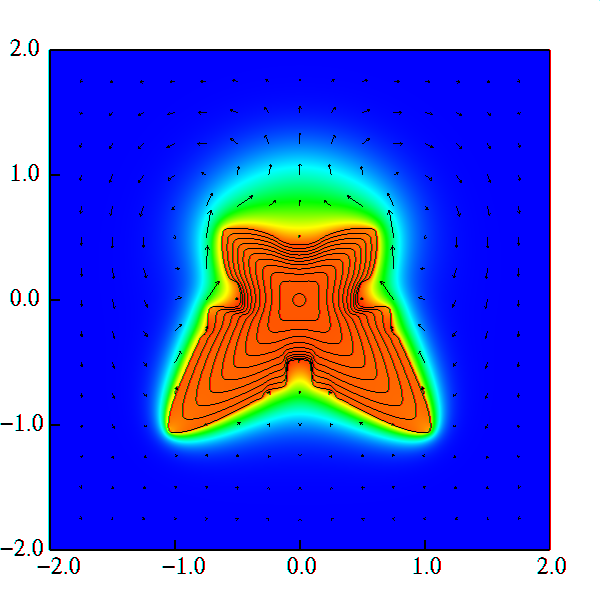} 
    \includegraphics[width=0.24\textwidth]{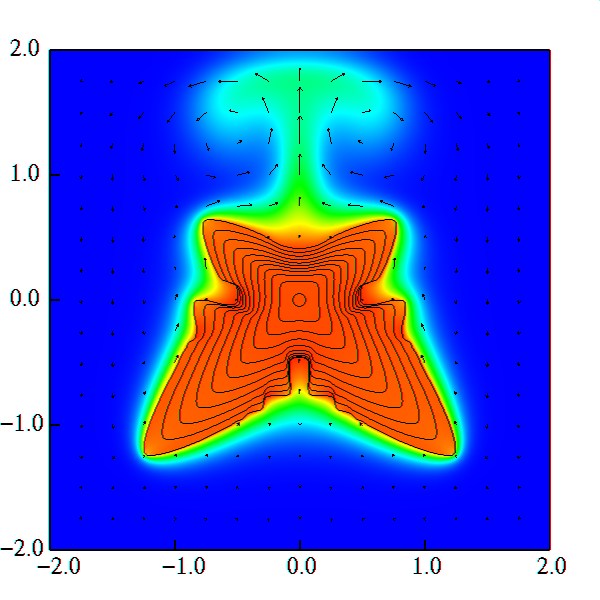} 
    \includegraphics[width=0.24\textwidth]{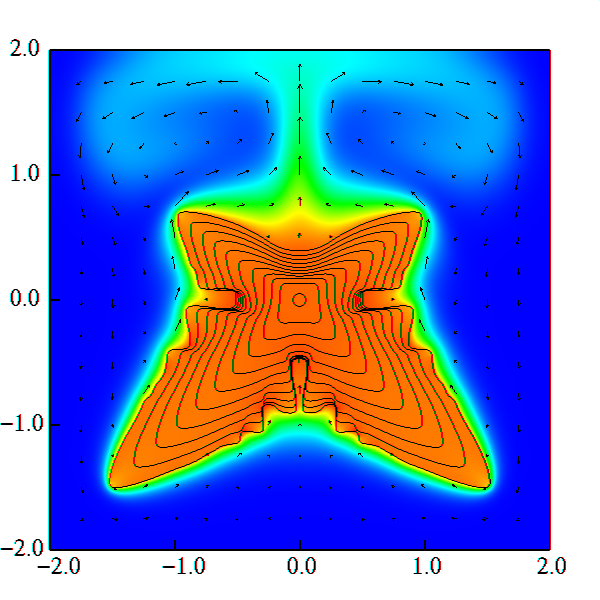} 
    \caption{Dendritic growth with buoyancy-driven flow. Interface shape, temperature field, and flow field at $t=0.1$ for $\beta = 0$, $10^3$, $2\times 10^3$, and $4\times 10^3$ from left to right.}
    \label{fig:buoyancy}
\end{figure}

In the final example, we consider the dendritic growth problem with buoyancy-driven flow.
The anisotropic surface tension is chosen in the rotated form $\varepsilon_C(\alpha) = 0.002(8/3\sin^4(2(\alpha-\pi/4)))$.
The no-slip boundary condition is applied to the fluid equation on all four boundaries.
We vary the thermal expansion coefficient so that the flow in the liquid phase is driven by buoyancy.
The gravitational acceleration is set to $g = 10$.
The reference temperature is set to the temperature of the surrounding fluid, $T_0 = -0.5$.
All other parameters are the same as in the previous example.
\Cref{fig:buoyancy} shows the results for increasing thermal expansion coefficients.
Near the solid-liquid interface, the released latent heat increases the fluid temperature and, as a result, causes fluid density changes and the buoyancy force.
Driven by buoyancy, the fluid transports heat from bottom to top, which produces an asymmetric temperature distribution in the $y$-direction.
The two upper branches are restrained from growing due to accumulated heat, while the two lower branches grow much faster since the heat flows away.

\section{Discussion}\label{sec:discuss}
This paper presents a numerical method for two representative moving interface problems based on a Cartesian grid-based boundary integral formulation coupled with an interface-evolution scheme.
The method combines the advantages of Cartesian grid solvers and boundary integral methods.

Elliptic and parabolic PDEs with irregular boundaries or interfaces are reformulated as well-conditioned boundary integral equations and solved with the KFBI method. The KFBI method uses a Cartesian grid-based solver to evaluate the integrals, which in most cases avoids the need for explicit analytical Green's functions and allows efficient use of fast PDE solvers such as FFTs and geometric multigrid methods. In the present work, however, the method is not completely kernel-free: a Green's function is still needed to impose artificial boundary conditions on the rectangular domain when the original problem is posed on an unbounded domain. Even so, the KFBI method avoids the evaluation of singular and nearly singular integrals, which are often difficult to treat accurately in quadrature-based boundary integral methods.

Interface evolution is computed in the $\theta-L$ formulation rather than the more common $x-y$ formulation. Because the interface is periodic, a Fourier pseudo-spectral method provides accurate spatial discretization. Moreover, the $\theta-L$ formulation enables a small-scale decomposition that removes the curvature-induced stiffness in both the Hele-Shaw flow and the Stefan problem. Combined with a semi-implicit scheme, this yields an efficient and stable FFT-based method for interface evolution in Fourier space.

While the current work primarily focuses on moving interface problems in two dimensions, certain improvements are needed to solve models in three dimensions, such as solidification problems and two-phase incompressible flows. First, elliptic PDEs with irregular boundaries and interfaces in three dimensions require a three-dimensional version of the KFBI method \cite{Ying2013}. Additionally, for parabolic PDEs, dimension-splitting techniques can be employed to accelerate computation \cite{Zhou2023}. Second, the $\theta-L$ approach is only applicable to evolving curves in two dimensions. To accurately track evolving surfaces in three dimensions, different numerical approaches, such as the front-tracking method or the level-set method, are needed. Finally, addressing the stiffness induced by mean curvature in three dimensions poses greater challenges than in two dimensions. Developing an efficient semi-implicit time-stepping scheme is crucial for ensuring the stable evolution of the interface.

\appendix
\section{Computation of derivative jumps}\label{sec:app1}
Suppose $\Gamma$ is parameterized as $\mathbf{X}(\alpha) = (x(\alpha), y(\alpha))$ where $\alpha$ is an arbitrary parameter.
Let $s$ be the arc-length parameter.
Suppose $\Gamma$ is sufficiently smooth, at least in the class $C^2$.
The unit outward normal is given by $\mathbf{n} = (y_{\alpha}/s_{\alpha}, -x_{\alpha}/s_{\alpha})$ where $s_{\alpha} = \sqrt{x_{\alpha}^2 + y_{\alpha}^2}$.
Given an interface problem with constant coefficients, we derive the jump values of the solution and its derivatives at the point $(x(\alpha), y(\alpha))\in\Gamma$.
\subsection{The modified Helmholtz equation}
Consider the interface problem of the modified Helmholtz equation
\begin{align}
\Delta u - c^2 u& = f, \quad\text{ in } \Omega^+\cup\Omega^-,\label{eqn:cdj-pde} \\ 
 [u] &= \Phi, \quad \text{ on } \Gamma, \label{eqn:cdj-ifc1} \\
[\partial_{\mathbf{n}} u] &= \Psi, \quad\text{ on }\Gamma. \label{eqn:cdj-ifc2}
\end{align}
The interface condition \eqref{eqn:cdj-ifc1} implies the zeroth-order jump value 
\begin{equation}
  [u] = \Phi.  
\end{equation}
By taking the derivative of both sides of the interface condition \eqref{eqn:cdj-ifc1} with respect to $\alpha$ and combining the interface condition \eqref{eqn:cdj-ifc2}, we have a $2\times 2$ linear system
\begin{align}
&x_{\alpha} [u_x] + y_{\alpha} [u_y] = \Phi_{\alpha}, \label{eqn:mh-jmp1-1}\\
&y_{\alpha}[u_x] - x_{\alpha}[u_y] = s_{\alpha}\Psi.\label{eqn:mh-jmp1-2}
\end{align}
Solving the linear system gives the values of $[u_x]$ and $[u_y]$.
By taking the derivative of both sides of \cref{eqn:mh-jmp1-1,eqn:mh-jmp1-2}, and using \cref{eqn:cdj-pde}, we have a $3\times 3$ linear system for $[u_{xx}]$, $[u_{yy}]$ and $[u_{xy}]$
\begin{align}
&(x_{\alpha})^2 [u_{xx}] + (y_{\alpha})^2 [u_{yy}] + 2x_{\alpha}y_{\alpha} [u_{xy}] = \Phi_{\alpha\alpha} -  x_{\alpha\alpha}[u_x] - y_{\alpha\alpha}[u_y],\\
&x_{\alpha}y_{\alpha}[u_{xx}] - x_{\alpha}y_{\alpha}[u_{yy}] + ((y_{\alpha})^2 - (x_{\alpha})^2)[u_{xy}] \nonumber\\
&= s_{\alpha\alpha}\Psi + s_{\alpha}\Psi_{\alpha} - y_{\alpha\alpha}[u_x] + x_{\alpha\alpha} [u_y],\\
&[u_{xx}] + [u_{yy}] = c^2 \Phi + [f].
\end{align}
After solving the linear system, the derivative jump values are obtained.

\subsection{The modified Stokes equation}
Consider the interface problem of the modified Stokes equation
\begin{align}
\Delta \mathbf{u} - c^2 \mathbf{u} - \nabla p & = \mathbf{f}, \quad\text{ in } \Omega^+\cup\Omega^-,\label{eqn:cdj-stokes} \\ 
\nabla \cdot \mathbf{u} & = 0, \quad\text{ in } \Omega^+\cup\Omega^-,\label{eqn:cdj-stokes-2} \\ 
 [\mathbf{u}] &= \mathbf{\Phi}, \quad \text{ on } \Gamma, \label{eqn:cdj-stokes-ifc1} \\
[T(\mathbf{u}, p)] &= \mathbf{\Psi}, \quad\text{ on }\Gamma. \label{eqn:cdj-stokes-ifc2}
\end{align}
where $\mathbf{u} = (u, v)^T$, $T(\mathbf{u},p) = -p\mathbf{n}+(\nabla\mathbf{u} + \nabla \mathbf{u}^T)\mathbf{n}$, $\mathbf{f} = (f_1,f_2)^T$, $\mathbf{\Phi} = (\Phi_1,\Phi_2)^T$, $\mathbf{\Psi} = (\Psi_1,\Psi_2)^T$.
The zeroth-order jump values are obtained from \cref{eqn:cdj-stokes-ifc1},
\begin{equation}
    [u] = \Phi_1, \quad [v] = \Phi_2.
\end{equation}
Taking the derivative of both sides of the interface condition \eqref{eqn:cdj-stokes-ifc1} with respect to $\alpha$ and using \cref{eqn:cdj-stokes-2,eqn:cdj-stokes-ifc2}, we obtain a $5\times 5$ system for the first-order jump values $[u_x]$, $[u_y]$, $[v_x]$, $[v_y]$ and $[p]$,
\begin{align}
    &x_{\alpha}[u_x] + y_{\alpha}[u_y] = \Phi_{1,\alpha},\label{eqn:ms-jmp1-1}\\
    &x_{\alpha}[v_x] + y_{\alpha}[v_y] = \Phi_{2,\alpha},\label{eqn:ms-jmp1-2}\\
    &2y_{\alpha}[u_x] - x_{\alpha}[u_y] - x_{\alpha}[v_x] - y_{\alpha}[p] = s_{\alpha}\Psi_{1},\label{eqn:ms-jmp1-3}\\
    &y_{\alpha}[u_y] + y_{\alpha}[v_x] - 2x_{\alpha}[v_y] + x_{\alpha}[p] = s_{\alpha}\Psi_{2},\label{eqn:ms-jmp1-4}\\
    &[u_x] + [v_y] = 0.
\end{align}
Taking the derivative of both sides of \cref{eqn:ms-jmp1-1,eqn:ms-jmp1-2,eqn:ms-jmp1-3,eqn:ms-jmp1-4} with respect to $\alpha$ and the derivatives of both sides of \cref{eqn:cdj-stokes-2} with respect to $x$ and $y$, and using \cref{eqn:cdj-stokes}, we obtain an $8\times 8$ system for the second-order jump values $[u_{xx}]$, $[u_{yy}]$, $[u_{xy}]$, $[v_{xx}]$, $[v_{yy}]$, $[v_{xy}]$, $[p_x]$ and $[p_y]$,
 \begin{align}
     &(x_{\alpha})^2[u_{xx}] + (y_{\alpha})^2[u_{yy}] + 2x_{\alpha}y_{\alpha}[u_{xy}] = r_1,\\
     &(x_{\alpha})^2[v_{xx}] + (y_{\alpha})^2[v_{yy}] + 2x_{\alpha}y_{\alpha}[v_{xy}] = r_2,\\
     &2x_{\alpha}y_{\alpha}[u_{xx}]-x_{\alpha}y_{\alpha}[u_{yy}]+(2(y_{\alpha})^2 -(x_{\alpha})^2)[u_{xy}] \nonumber\\ &-(x_{\alpha})^2[v_{xx}]-x_{\alpha}y_{\alpha}[v_{xy}]-x_{\alpha}y_{\alpha}[p_x]-(y_{\alpha})^2[p_y] = r_3, \\
     &(y_{\alpha})^2[u_{yy}] + x_{\alpha}y_{\alpha}[u_{xy}]+x_{\alpha}y_{\alpha}[v_{xx}]-2x_{\alpha}y_{\alpha}[v_{yy}]\nonumber\\
     &+((y_{\alpha})^2 - 2(x_{\alpha})^2)[v_{xy}]+(x_{\alpha})^2[p_x]+x_{\alpha}y_{\alpha}[p_y] = r_4,\\
     &[u_{xx}] + [u_{yy}] - [p_x] = c^2 \Phi_1 + [f_1],\\
     &[v_{xx}] + [v_{yy}] - [p_y] = c^2 \Phi_2 + [f_2],\\
     &[u_{xx}] + [v_{xy}] = 0,\\
     &[u_{xy}] + [v_{yy}] = 0.
 \end{align}
 where $r_i, i = 1,2,\cdots,4$ are given by
 \begin{align}
 &r_1 = \Phi_{1,\alpha\alpha} - x_{\alpha\alpha}[u_x] - y_{\alpha\alpha}[u_y],\\
 &r_2 = \Phi_{2,\alpha\alpha} - x_{\alpha\alpha}[v_x] - y_{\alpha\alpha}[v_y],\\
 &r_3 = s_{\alpha\alpha}\Psi_{1}+s_{\alpha}\Psi_{1,\alpha}-2y_{\alpha\alpha}[u_x]+x_{\alpha\alpha}[u_y]+x_{\alpha\alpha}[v_x]+y_{\alpha\alpha}[p],\\
 &r_4 = s_{\alpha\alpha}\Psi_{2}+s_{\alpha}\Psi_{2,\alpha} - y_{\alpha\alpha}[u_y]-y_{\alpha\alpha}[v_x]+2x_{\alpha\alpha}[v_y]-x_{\alpha\alpha}[p].
 \end{align}
By solving the three linear systems, the derivative jump values of $u$, $v$, and $p$ can be obtained.

\section*{Acknowledgement}
W.~Ying was supported by the National Natural Science Foundation of China, Division of Mathematical Sciences (Project No.~12471342), and by the Fundamental Research Funds for the Central Universities of China.
S.~Li was partially supported by the National Science Foundation (NSF), United States of America under grant DMS-2309798.

\bibliographystyle{elsarticle-num} 
\bibliography{mybibfile}





\end{document}